\documentclass[review]{elsarticle}

\usepackage{amsmath,amsthm,amsfonts,amssymb,amscd,bigints}
\usepackage{subfigure}
\usepackage{graphicx,bm,color}
\usepackage{algorithm}
\usepackage{algpseudocode}
\usepackage{stmaryrd}
\usepackage{listings}
\usepackage{lineno}
\usepackage{mathtools}
\usepackage[version=3]{mhchem}
\usepackage[colorlinks,allcolors=blue]{hyperref}

\usepackage{lineno}

\makeatletter
\def\ps@pprintTitle{%
 \let\@oddhead\@empty
 \let\@evenhead\@empty
 \def\@oddfoot{}%
 \let\@evenfoot\@oddfoot}
\makeatother

\definecolor{gray}{gray}{0.6}

\def\u{{\bm u}}

\def\g{{\bm g}}

\def\grad{\nabla}
\def\div{\nabla \cdot}

\begin{document}

\title{A Domain Decomposition Approach for Local Mesh Refinement in Space and Time}

\author[csm]{Gurpreet~Singh\corref{cor1}}
\ead{gurpreet@utexas.edu}
\author[csm]{Mary~F.~Wheeler}
\ead{mfw@ices.utexas.edu}

\cortext[cor1]{Corresponding author}
\address[csm]{Center for Subsurface Modeling, The University of Texas at Austin, Austin, TX 78712}

\begin{abstract}
Reservoir simulations for complex multiphase flow and transport problems often suffer from non-linear solver convergence issues. These manifest in the form of restrictively small time-step sizes even while using unconditionally stable fully implicit schemes. These problems are further compounded when a local mesh refinement is used to accurately represent reservoir parameters available such as permeability, porosity, etc., at multiple spatial scales. We discuss a domain decomposition approach that allows different time-step sizes and mesh refinements in different subdomains \cite{SinghEVST} of the reservoir that circumvents these issues without compromising computational efficiency and prediction accuracy. This approach extends the well-known methodology of local mesh refinement in space \cite{WheelerEV} to time. Our numerical experiments indicate that non-linear solvers fail to converge, to the desired tolerance, due to large non-linear residuals in a smaller subdomain. We exploit this feature to identify subdomains where smaller time-step sizes are necessary while using large time-step sizes in the rest of the reservoir domain. The three key components of our approach are: (1) a space-time, enhanced velocity, domain decomposition approach that allows different mesh refinements and time-step sizes in different subdomains while preserving local mass conservation, (2) a residual based error estimator to identify or mark regions (or subdomains) that pose non-linear convergence issues, and (3) a fully coupled monolithic solver is also presented that solves the coarse and fine subdomain problems, both in space and time, simultaneously. This solution scheme is fully implicit and is therefore unconditionally stable. The results indicate that using large time-step sizes for the entire reservoir domain poses serious non-linear solver convergence issues. Although using a smaller time step size for the entire domain reduces the convergence issues, it also results in substantial computational overheads. The proposed space-time domain decomposition approach, with smaller time-step sizes in a subdomain and large time-step sizes everywhere else, circumvents the non-linear convergence issue without adding computational costs. Additionally, a space-time monolithic solver renders a massively parallel, time concurrent framework for solving flow and transport problems in subsurface porous media. Since the proposed approach is similar to the widely used finite difference scheme, it can be easily integrated in any existing legacy reservoir simulator.
\end{abstract}

\begin{keyword}
dynamic mesh refinement \sep space-time domain decomposition \sep mixed finite element  \sep enhanced velocity \sep monolithic system \sep fully-implicit
\end{keyword}

\maketitle



\newcommand{\bs}[1]{\boldsymbol{#1}}


\section{Introduction}

In this section we outline a framework for space-time adaptivity for a two-phase flow example followed by an abstraction that allows extensions to be general multiphase flow and reactive transport problems in porous media. The purpose of this section is to also present a strong motivation for the development of our proposed space-time adaptive mesh refinement approach.
\subsection{Tracking Features}\label{subsec:features}
\begin{figure}[H]
\begin{center}
\includegraphics[width=14cm]{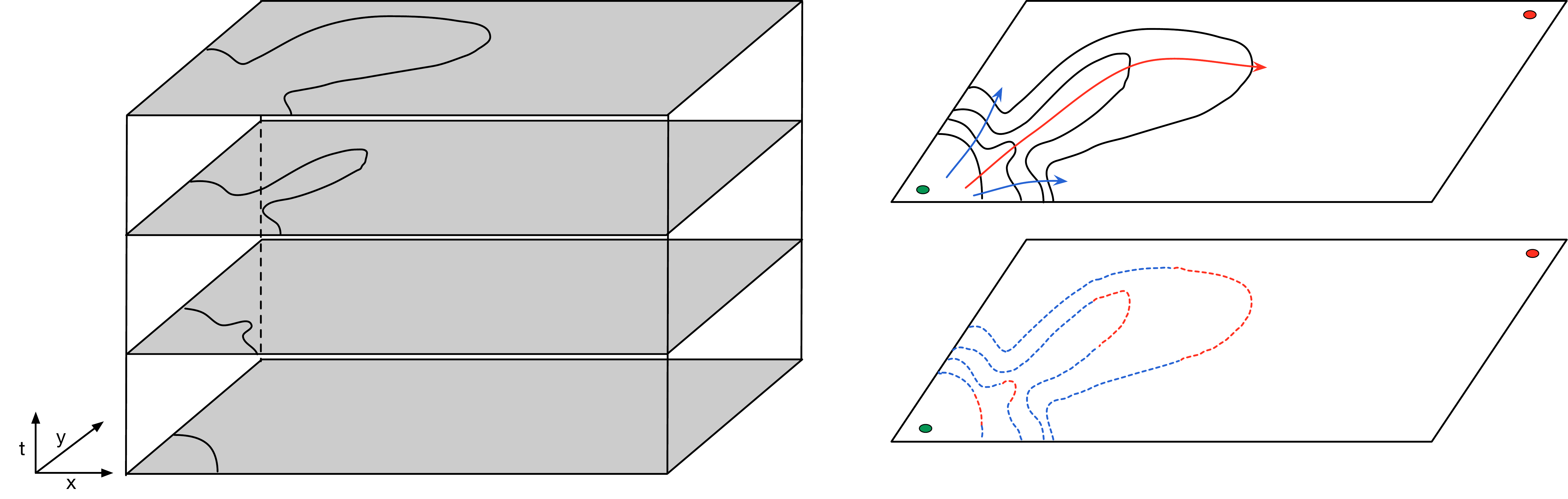}
\caption{Example of saturation front evolution in space and time over a reservoir domain.}
\label{fig:example}
\end{center}
\end{figure}

Let us consider a two-phase (oil-water) flow problem as shown in Figure \ref{fig:example} with injection and production wells located at the bottom left and top right corners of the reservoir domain. Here, we assume two spatial dimensions, indicated by x and y, in the horizontal plane and one temporal dimension in the vertical direction. Figure \ref{fig:example} (left) shows a hypothetical evolution of a saturation front ($S^{*}_{w}$ = constant) in space and time. Although, the approach presented in this work has been generalized to three spatial dimensions, we restrict ourselves to just two for simplicity. Figure \ref{fig:example} (top right) shows the iso-saturation maps for the different times in Figure \ref{fig:example} (left). With this description, we now examine the quantity of interest or the saturation using a delta operator to identify change in space ($\Delta_{s} $) and time ($\Delta_{t}$). We are well aware that the magnitude of $\Delta_{s} S_{w}$ is highest in the vicinity of the saturation front at a given time. Furthermore, from the iso-saturation map in Figure \ref{fig:example} (top right) we also note that the magnitude of ($\Delta_{t} S_{w}$) is highest along the red line. Let us now consider a simple task of dividing the reservoir into non-overlapping subdomain using these changes in saturation in space and time. A simple binary classification in terms of slow and fast changes of saturation allows us to identify three regions in Figure \ref{fig:example} (top right): (1) The red line along which the magnitudes of both $\Delta_{t} S_{w}$ and $\Delta_{s} S_{w}$ are high, (2) the blue lines along which the magnitude of $\Delta_{t} S_{w}$ is small and $\Delta_{s} S_{w}$ is large, and (3) the remaining reservoir domain where both $\Delta_{t} S_{w}$ and $\Delta_{s} S_{w}$ are small. Figure \ref{fig:example} (right, bottom) shows a further classification of the iso-saturation map into these three subdomain descriptions. We now consider the following questions:
\begin{enumerate}
\item{How can we preserve solution accuracy with the saturation front evolving at multiple time-scales in different parts of the reservoir domain?}
\item{What are the computational bottlenecks and where are they located in the space-time reservoir domain?}
\item{Can such a classification be used to balance the computational loads and improve efficiency?}
\end{enumerate}

An accurate representation of the saturation front at a given time requires a local spatial mesh refinement with the refinement factor depending on the desired accuracy or availability of reservoir petrophysical properties (permeability, porosity etc.) at different spatial scales. For example, data from well logs and geological models is often available at a resolution of a few feet and in order to incorporate such fine scale information local, spatial mesh refinement is often necessary to preserve solution accuracy. However, to promote computational efficiency one can use adaptive upscaling approaches \cite{Yerlan17,singhatce17,singhcg18} wherein a local spatial mesh refinement is used only along the saturation front while using upscaled properties away from it without loss in solution accuracy. In doing so, we are tracking the saturation front as it evolves spatially achieving the desired accuracy. However, we have not sufficiently addressed the multiple time-scales issue in the first question. In a heterogeneous reservoir domain, the saturation front moves at different velocities in different reservoir subdomains due local permeability differences. Analogous to local spatial mesh refinement along the saturation front for accuracy it also makes sense if we could take smaller time-step sizes (local temporal mesh refinement) along the saturation front where the changes in time are large (large $\Delta_{t} S_{w}$). This is shown as red dotted line on the iso-saturation map in Figure \ref{fig:example} (bottom right). Using the aforementioned classification of the reservoir domain based upon changes in saturation in space and time we now obtain three local mesh refinements (1) fine in space and time where both $\Delta_{s} S_{w}$ and $\Delta_{t} S_{w}$ are high (red dotted lines), (2) fine in space and coarse in time where $\Delta_{s} S_{w}$ is high but $\Delta_{t} S_{w}$ is low (blue dotted lines), and (3) coarse in space and time where both $\Delta_{s} S_{w}$ and $\Delta_{t} S_{w}$ are low (remaining reservoir domain). With this, we have roughly outlined a strategy towards addressing the solution accuracy as the problem evolves in space and time. 

To answer the second question, let us first define a weak metric for computational cost to identify problem areas or bottlenecks. A large array of subsurface porous media problems are non-linear due to property descriptions such as relative permeability, density etc., to name a new, as non-linear functions of the primary variables (pressure, saturation, concentration etc.). Numerical reservoir simulations commonly use either a Newton's method or it's variants \cite{kelley1995iterative} to linearize the non-linear system of algebraic equations resulting from spatial (FDM, FEM, FVM; Finite Difference, Finite Element Methods, and Finite Volume Methods)  and temporal discretization (backward Euler, forward Euler or Crank-Nicolson) of the partial differential equations (PDE) in a given model formulation. A linear solver is then used to solve this resulting system of now linear algebraic equations. For the sake of defining a weak metric for computational cost, let us assume that for a given system (after Newton linearization) the computational cost of solving the linear system (linear solver cost) remains the same for a fixed degrees of freedom (DOF). In other words, the linear solver computational cost is linearly related to the degrees of freedom. We can now characterize the computational load by the number of Newton (or other non-linear solver) iterations as the weak metric. For a linear system such as incompressible, single phase flow and tracer transport in the absence of velocity dependent dispersion, the computational load is already balanced in terms of this metric. A linear system always converges in one Newton iteration with any spatial and temporal mesh refinement. Here, a local spatial and temporal mesh refinement only contributes to the solution accuracy and does not create computational bottlenecks. 

However, for a non-linear problem the solution accuracy and computational load are strongly tied to each other. For example, in a two-phase, compressible flow problem this non-linearity manifests itself at the saturation front where the changes in saturation are large and hence the resulting non-linear, relative permeability changes are large. In terms of our weak computational metric (Newton iterations), this results in additional Newton iterations leading to increased computational costs. However, it is interesting to note that away from the saturation front (spatially) where changes in saturation (temporally) are small, the computational load must also remain small. Furthermore, part of the saturation front that evolves slowly in time the previous time solution remains close to the next time solution (Figure \ref{fig:example}, bottom right, dotted blue part of the saturation iso-map) and therefore an analogous argument of low computational load applies here as well. 

\begin{figure}[H]
\begin{center}
\includegraphics[width=7.5cm,trim=2cm 3cm 2cm 2cm, clip]{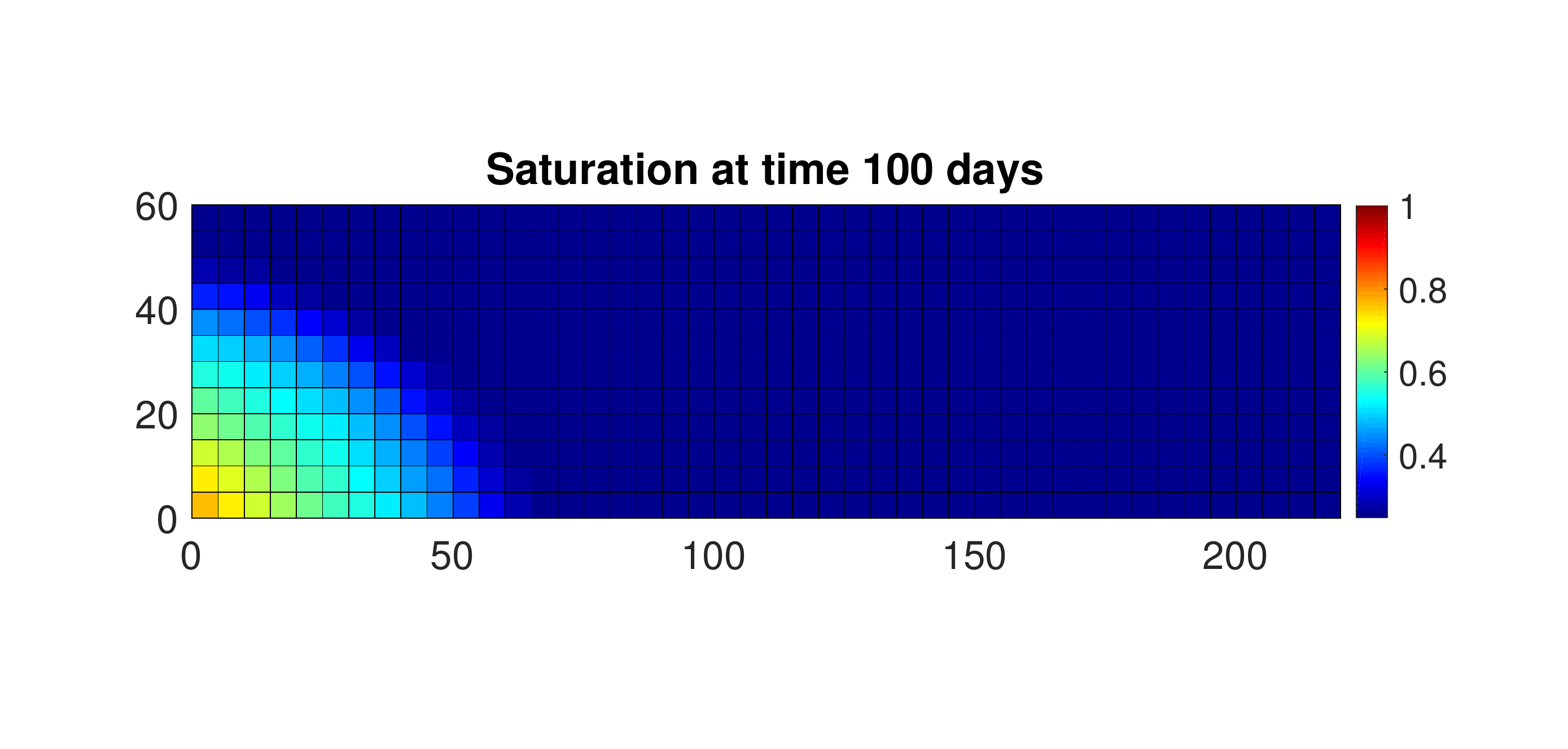}
\includegraphics[width=7.5cm,trim=2cm 3cm 2cm 2cm, clip]{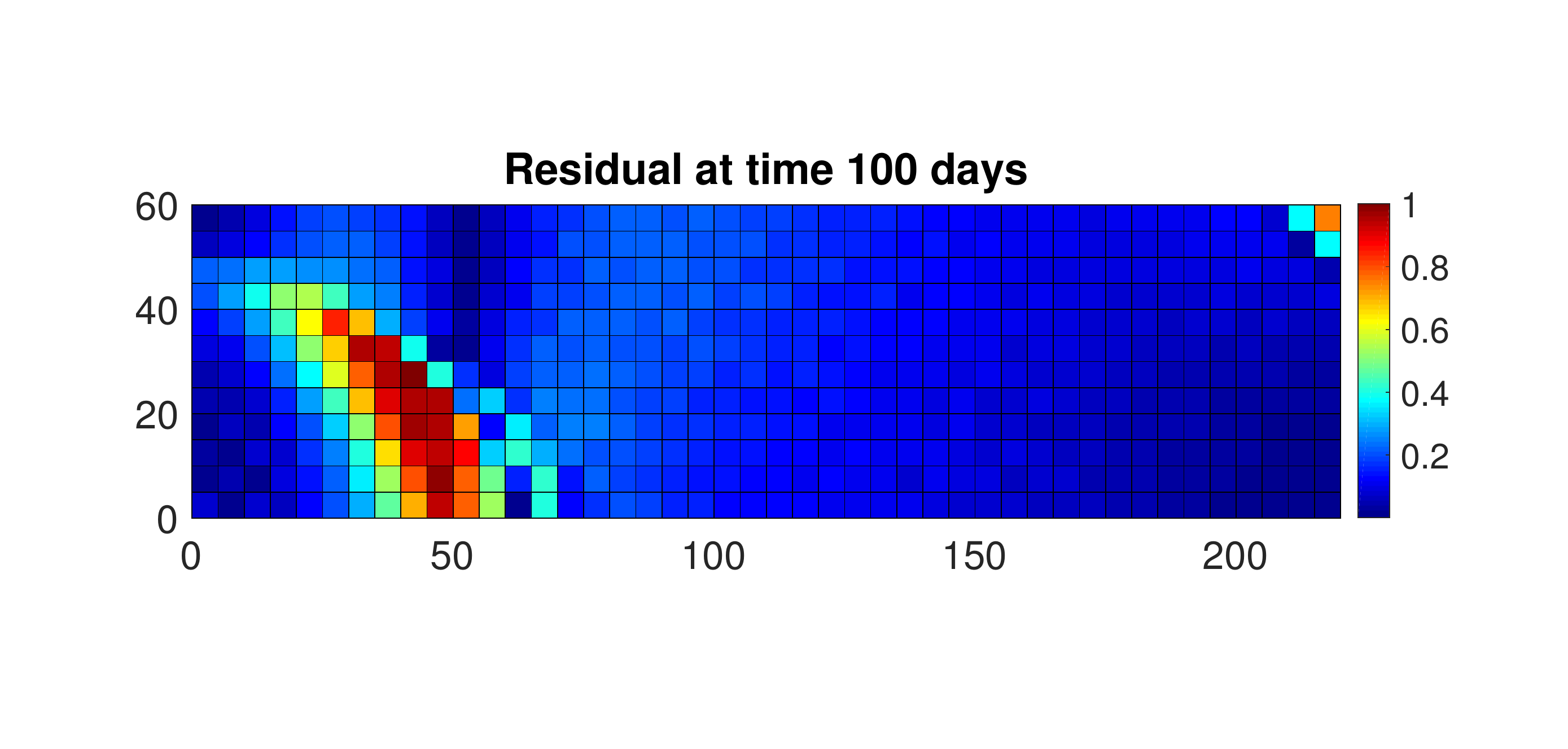}\\
\includegraphics[width=7.5cm,trim=2cm 3cm 2cm 2cm, clip]{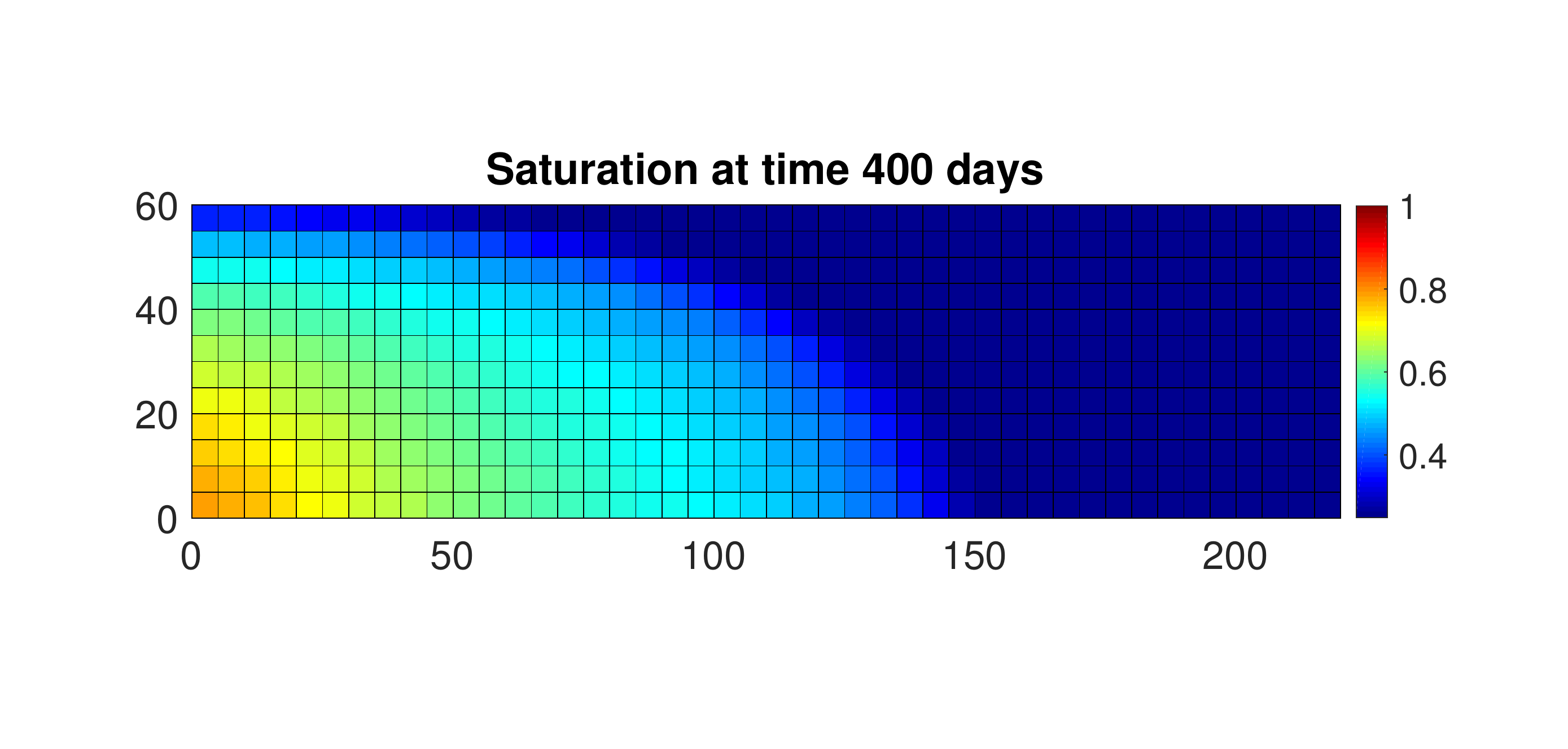}
\includegraphics[width=7.5cm,trim=2cm 3cm 2cm 2cm, clip]{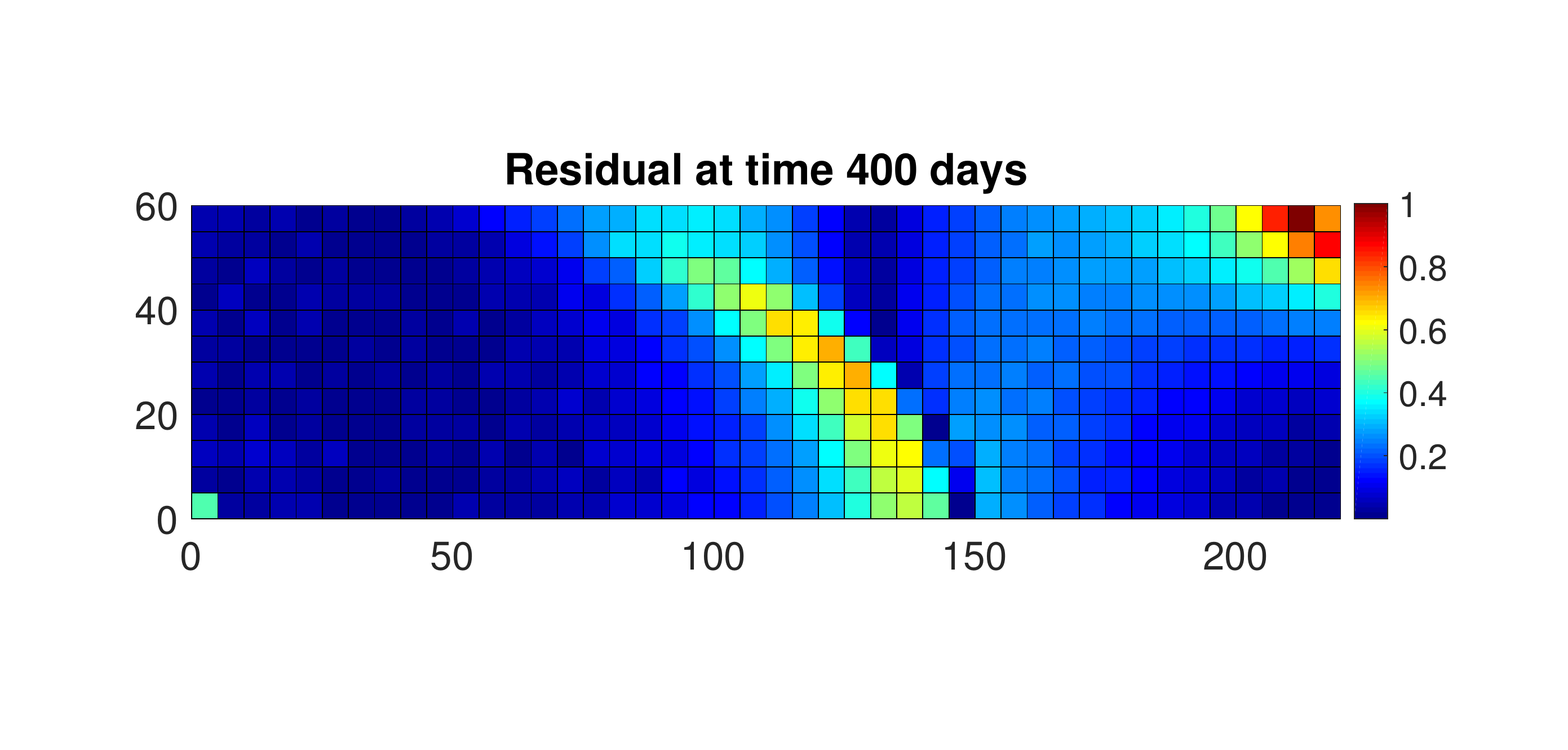}
\caption{Saturation distribution and normalized, non-linear residuals at 100 and 400 days.}
\label{fig:resid}
\end{center}
\end{figure}

A spatio-temporal distribution of the non-linear normalized residuals for a two-phase flow problem prior to the non-linear system convergence to a desired tolerance identifies this region of high computational load. Figure \ref{fig:resid} shows the saturation distributions at 100 and 400 days and the corresponding distributions of non-linear residuals for a two-phase, slightly compressible flow problem in a homogeneous porous medium. Here, high non-linear residuals are located at the saturation front and near the injection and production wells located at the diagonally opposite corners of the reservoir domain. We now make a conjecture that the number of Newton iterations (the weak metric) for a system to converge to a desired tolerance is dominated by a narrow region in space and time where the non-linearity is high. Please note that the weak computational metric defined above assumes a linear dependency of computational cost on the degrees of freedom. This assumption allows us to focus on the computational bottlenecks due to non-linearities alone and can be easily relaxed by using non-linear solver tolerance to specify linear solver tolerance using forcing function approaches (\cite{eisenstat})  and specialized preconditioners (\cite{lacroix2001decoupling,lacroix2003iterative,singhapprox,Klie96two-stagepreconditions}).

To answer the third question, we will rely upon a recently developed space-time domain decomposition approach \cite{SinghEVST} as a vehicle that uses the above information to promote overall computational efficiency. This approach is locally mass conservative and allows for different mesh refinements and time-step sizes in different subdomains of the reservoir. A cheap explicit error estimator based upon non-linear residuals is then used to identify regions for local refinements in space and time in a small subdomain of the reservoir while keeping the mesh coarse in the remaining reservoir domain. A number of error estimators \cite{Vohralik13} and references therein) are available in the literature that distinguish between spatial, temporal, and solver (linear and non-linear) errors to guide dynamic mesh refinement strategies. This computational load balancing strategy results in substantial speedups (up to a factor of 40 compared to conventional approaches) in a prototype implementation of this framework on a serial machine. This framework and our space-time domain decomposition approach are inherently flexible and can be adapted to different parallel compute architectures (shared or distributed). In this work, we restrict ourselves to this serial prototype for computational benchmarking and ease of describing the key elements of this framework. A parallel implementation with problem specific, explicit error estimator design is being actively pursued and will be presented in future works.

\subsection{Tracking Non-linearities}\label{subsec:nonlin}

So far, we discussed the adaptive space-time framework in the light of a two-phase flow problem. However, this framework can be generalized and extended to handle general multiphase flow and reactive transport problems such as black oil and compositional flow for water, gas and chemical flooding processes. Here, we briefly describe an abstraction of the above framework in order to outline feasible paths to future extensions. Invoking the arguments of delta change in a quantity of interest ($Q$) and the binary classification we can decompose the space-time reservoir domain into four subdomains as shown in Table 1. We can now apply the above framework to a wide range of problems by defining this quantity of interest. For the two-phase problem we only used the identifiers 1, 2, and 4 since they allow us to adequately balance the computational load during simulation runtimes. Identifier 3 is useful for problems involving reactive chemistry (kinetic or equilibrium) and other local equilibrium calculations. For example, during gas flooding in a heterogeneous reservoir miscibility can vary spatially from first (or multi-contact) miscibility to immiscible states. Here, first contact and immiscible fronts are in the vicinity of a concentration front and hence can be identified by using a combination of delta change ($\Delta_{s,t}$) in concentrations and explicit estimators (normalized residuals). However, multi-contact miscibility regions are often characterized by large concentration changes ($\Delta_{t}$) in time with small concentrations changes in space ($\Delta_{s,t}$) indicated by the identifier 3 depending on the nature of the local phase behavior calculations for a given hydrocarbon fluid composition. The same is true for surface (or heterogeneous) reactions wherein the changes in time might be large behind the concentration front due to high non-linear reaction rates. Table 1 provides a simple two-scale classification (coarse and fine) of the space-time domain towards distributing computational loads based upon domain knowledge of the problem at hand. 

As for the quantify of interest ($Q$), the definition varies and is subjective to the problem at hand. For example, for a reactive transport problem with multiple kinetic reactions and rate constants varying orders of magnitude, it makes sense to take smaller time-steps for fast reactions and larger time-steps for slower reactions without compromising accuracy. Therefore, the quantify of interest here is the reaction rate itself. A naive application of explicit estimators combined with delta change in concentration is not sufficient to adequately balance the computational load here. In fact, the explicit estimators (non-linear residuals at each Newton iteration) lump the errors due to flow, transport, and reactions into a single quantity making it difficult to resolve which amongst the three is dominant. Therefore, knowledge of physical processes (domain knowledge) is necessary to construct problem specific error estimators.  This abstract framework now allows us to bring common-sense domain knowledge (understanding of physical processes) towards improving overall computational efficiency of numerical simulations.

With this outline for the space-time adaptivity framework, we first describe the domain decomposition approach that allows us to take different spatial mesh refinements and time-step sizes in different subdomains of a large reservoir. Here we briefly discuss the variational (weak) formulation of a system of PDEs representing the model problem. This is followed by a description of a time-concurrent, monolithic, space-time solver for the resulting system of non-linear algebraic equations resulting from the space-time discretization. In this section, we differentiate between time-marching nature of conventional solution schemes and the time-concurrent aspect of our proposed solution algorithm. We then present a set of numerical results that demonstrate speedups by factors of 25 to 40 times for our space-time adaptive framework over conventional solution schemes (serial computations for both) in solving a multiphase flow and transport problem. Finally, we present a summary of the presented work with conclusions and future outlook towards development of a parallel, space-time adaptive framework and associated development for error estimators.

\begin{table}
  \begin{center}
    \caption*{Table 1: Space-time domain decomposition based upon change in a quantity of interest $Q$}
    \begin{tabular}{|c|c|c|l|} 
     \hline
      \textbf{Identifier} & \textbf{Change in space}  & \textbf{Change in time} & \textbf{Local Mesh}\\
      $\alpha$ & $\Delta_{s} Q$ & $\Delta_{t}Q$ & $h$\\ 
      \hline
      1 & large & large & fine in space and time\\
      2 & large & small & fine in space, coarse in time\\
      3 & small & large & coarse in space, fine in time\\
      4 & small & small & coarse in space and time\\
       \hline
    \end{tabular}
  \end{center}
\end{table}

\section{Space-Time Enhanced Velocity Domain Decomposition}\label{sec:evst}
In this section, we briefly describe the space-time enhanced velocity domain decomposition approach applied to a slightly compressible, two-phase (oil-water) problem in porous media. For further details regarding the fully discrete variational formulation, solution spaces, and quadrature rules the reader is referred to \cite{SinghEVST}. Here our primary focus is on constructing a non-linear algebraic system by applying the aforementioned domain decomposition on the PDEs associated with the model formulation. To this effect, we first describe the two-phase flow model formulation, followed by a short description of the weak variational form, and evaluation of appropriate integrals to obtain a non-linear algebraic system of equations. Please note that the choice of solution spaces and quadrature rules inherent to this domain decomposition scheme results in a discretization scheme that is very similar to the widely used finite difference method. Therefore, this scheme can be incorporated in any FDM based legacy reservoir simulator with minimal changes to the existing code framework. 

\subsection{Two-phase flow formulation}\label{subsec:formulation}
We present the well-known immiscible, two-phase, slightly compressible flow in porous medium model formulation with oil and water phase mass conservation and constitutive equations along with the boundary and initial conditions. The mass conservation equation for phase $\alpha$ is given by,
\begin{equation}
\frac{\partial\left(\phi\rho_{\alpha}s_{\alpha}\right)}{\partial t} + \div \u_{\alpha} = q_{\alpha} \text{~in~} \Omega \times J,
\label{eqn:2phcon}
\end{equation}
where $\phi$ and $K$ have their usual meanings as described before, and $\rho_{\alpha}$, $s_{\alpha}$, $u_{\alpha}$ and $q_{\alpha}$ are density, saturation, velocity and source/sink term, respectively of phase $\alpha$. The constitutive equation for the corresponding phase $\alpha$ is given by Darcy's law as,
\begin{equation}
\u_{\alpha} = -K\rho_{\alpha}\frac{k_{r\alpha}}{\mu_{\alpha}}\left(\grad p_\alpha-\rho_{\alpha}\g\right) \text{~in~} \Omega \times J
\end{equation}
Further, $k_{r\alpha}$, $\mu_{\alpha}$ and $p_{\alpha}$ are the relative permeability, viscosity and pressure of phase $\alpha$. 
Although not restrictive, for the sake of simplicity we assume no flow boundary conditions.
\begin{eqnarray}
\u_{\alpha} \cdot \bs{\nu} = 0 \text{~on~} \partial \Omega \times J\\
p_{\alpha} = p_{\alpha}^{0}, \quad s_{\alpha} = s_{\alpha}^{0},  \text{~at~} \Omega \times \{0\}
\end{eqnarray}
Here, $p^{0}_{\alpha}$, $s_{\alpha}^{0}$ are the initial conditions for pressure and saturation of phase $\alpha$.
Furthermore, the phase saturations $s_{\alpha}$ obey the constraint,
\begin{equation}
\sum_{\alpha}s_{\alpha} = 1.
\label{eqn:2phsat}
\end{equation}
We assume capillary pressure and relative permeabilities to be continuous and monotonic functions of phase saturations, 
\begin{equation}
p_{c}=f(S_{w}) = p_{o}-p_{w},
\label{eqn:2phcap}
\end{equation}
\begin{equation}
k_{r\alpha} = k_{r\alpha}(s_{\alpha}).
\end{equation}
The oil and water phase are assumed to slightly compressible with phase densities evaluated using,
\begin{equation}
\rho_{\alpha} = \rho_{\alpha,ref}\exp\left[c_{f\alpha}(p_{\alpha}-p_{\alpha,ref})\right].
\end{equation}
Here, $c_{f\alpha}$ is the compressibility and $\rho_{\alpha,ref}$ is the density of phase $\alpha$ at the reference pressure $p_{\alpha,ref}$.

\subsection{Fully discrete formulation (or non-linear algebraic system)} 
In this subsection, we briefly touch upon the weak variational and fully discrete forms to obtain a non-linear system of algebraic equations in one spatial and one temporal dimension. For complete details regarding the choice of spaces, quadrature rules, and the expanded mixed variational formulation leading to the space-time enhanced velocity domain decomposition approach the reader is referred to \cite{SinghEVST}. The expanded mixed variational form of Eqns. \eqref{eqn:2phcon} thru \eqref{eqn:2phsat} is: find $\bs{u}_{\alpha,h}^{t} \in \bs{V}_{h}^{t,*}$,  $\tilde{\bs{u}}_{\alpha,h}^{t} \in \bs{V}_{h}^{t,*}$, $s_{w,h}^{t} \in W_{h}^{t}$, and $p_{o,h}^{t} \in W_{h}^{t}$ such that,

\begin{equation}
\left(\frac{\partial}{\partial t} \phi\left(\rho_{w}s^{t}_{w,h} + \rho_{o}(1-s^{t}_{w,h})\right),w\right) + \left(\div \left(\bs{u}^{t}_{w,h}+\bs{u}^{t}_{o,h}\right),w\right) = \left(q_{w}+q_{o},w\right)
\label{eqn:totcon}
\end{equation}

\begin{equation}
\left(\frac{\partial}{\partial t} \left(\phi \rho_{w}s^{t}_{w,h} \right),w\right) + \left(\div \bs{u}^{t}_{w,h},w\right) = \left(q_{w},w\right),
\label{eqn:watcon}
\end{equation}

\begin{equation}
\left(K^{-1}\tilde{\u}^{t}_{o,h},\bs{v}\right) - \left(p^{t}_{o,h},\div \bs{v}\right) =0,
\label{eqn:oildar}
\end{equation}
\begin{equation}
\begin{aligned}
\left(K^{-1}\tilde{\u}^{t}_{w,h},\bs{v}\right) - \left(p^{t}_{o,h},\div \bs{v}\right) = -\left(p_{c},\div \bs{v}\right) ,
\end{aligned}
\label{eqn:watdar}
\end{equation}
\begin{equation}
\begin{aligned}
\left( \u^{t}_{\alpha,h},\bs{v}\right) = \left(\lambda_{\alpha} \tilde{\u}^{t}_{\alpha,h},\bs{v}\right),
\end{aligned}
\label{eqn:expand}
\end{equation}
with $w\in W$ and $\bs{v}\in \bs{V}$. Please note that $s_{o}$ and $p_{w}$ are eliminated in the above formulation using the algebraic constraints Eqns. \eqref{eqn:2phsat} and \eqref{eqn:2phcap}, respectively in favor of $s_{w}$ and $p_{o}$. Further, $\lambda_{\alpha}$ is defined as the mobility of phase $\alpha$ as,
\begin{equation}
\lambda_{\alpha}= \frac{k_{r\alpha}\rho_{\alpha}}{\mu_{\alpha}},
\end{equation}
An expanded mixed formulation \cite{weiser85, malgo02}, with additional auxiliary phase fluxes $\tilde{\bs{u}}_{\alpha}$, is used to avoid inverting zero phase relative permeabilities ($k_{r\alpha}$). The solution can then be written as,
\begin{equation}
p_{o} = \sum_{m=1}^{q}\sum_{i = 1}^{r} P_{i}^{m} w_{i}^{m}, \quad \bs{u}_{\alpha} = \sum_{m=1}^{q}\sum_{i = 1}^{r+1} U_{\alpha,i+\frac{1}{2}}^{m} \varphi_{i+\frac{1}{2}}^{m},
\end{equation}
\begin{equation}
s_{w} = \sum_{m=1}^{q}\sum_{i = 1}^{r} S_{w,i}^{m} w_{i}^{m}, \text{ and} \quad \bs{\tilde{u}}_{\alpha} = \sum_{m=1}^{q}\sum_{i = 1}^{r+1} \tilde{U}_{\alpha,i+\frac{1}{2}}^{m} \varphi_{i+\frac{1}{2}}^{m},
\end{equation}
We will now construct an algebraic system of equations by testing the variational forms of the discrete constitutive and conservation equations with $w_{j}^{n}$ and $\varphi_{j+\frac{1}{2}}^{n}$, respectively. Here, we evaluate most of the integral terms in the variational problem on a matching grid, bifurcating to the non-matching grid only for the integrals where the non-matching, space-time interface enters the evaluation. Again for simplicity, let us assume that the coarse ($\delta t_{c}$) and fine ($\delta t_{f}$) time step sizes such that the ratio $\delta t_{c}/\delta t_{f} = 3$ indicating that the fine time subdomain is three times refined with respect to the coarse subdomain. on  For the first terms in the constitutive Eqns. \eqref{eqn:oildar} and \eqref{eqn:watdar} we have,

\begin{equation}
\begin{aligned}
\left( K^{-1}\tilde{\bs{u}}_{\alpha},\varphi_{j+\frac{1}{2}}^{n}\right)_{\Omega\times J} & \approx  \left( K^{-1}\sum_{m=1}^{q} \sum_{i=1}^{r+1} \tilde{U}_{\alpha,i+\frac{1}{2}}^{m}\varphi_{i+\frac{1}{2}}^{m}, \varphi_{j+\frac{1}{2}}^{n} \right)_{TM}\\
& = \frac{1}{2 \left| e_{j+\frac{1}{2}}^{n} \right|} \left( \frac{h_{j}}{K_{j}} + \frac{h_{j+1}}{K_{j+1}} \right) U_{\alpha,j+\frac{1}{2}}^{n}
\end{aligned}
\end{equation}

\begin{equation}
h_{j} = x_{j+\frac{1}{2}} - x_{j-\frac{1}{2}}
\end{equation}
Here, $h_{j} = x_{j+\frac{1}{2}} - x_{j-\frac{1}{2}}$ and $\alpha$ is either oil ($o$) or water ($w$) phase. The second terms in Eqns. \eqref{eqn:oildar} and \eqref{eqn:watdar} can be expanded as,
\begin{equation}
\begin{aligned}
\left(p_{o},\nabla \cdot \varphi_{j+\frac{1}{2}}^{n}\right)_{\Omega\times J} & = \left(\sum_{m=1}^{q}\sum_{i=1}^{r} P_{o,i}^{m} w_{i}^{m},\nabla \cdot \varphi_{j+\frac{1}{2}}^{n}\right)_{\Omega\times J} \\
& = P_{o,j}^{n} - P_{o,j+1}^{n}.
\end{aligned}
\end{equation}
Let us now consider, for a given $j_{0}\in \mathbb{I}$, a non-matching grid with fine a domain at ${\left(j_{0}+\frac{1}{2}\right)}^{-}$ and a coarse domain at ${\left(j_{0}+\frac{1}{2}\right)}^{+}$,
\begin{equation}
\begin{aligned}
\left(p_{o},\nabla \cdot \varphi_{j_{0}+\frac{1}{2}}^{n-\frac{1}{3}}\right)_{\Omega\times J} & = \left(\sum_{m=1}^{q}\sum_{i=1}^{r} P_{o,i}^{m} w_{i}^{m},\nabla \cdot \varphi_{j_{0}+\frac{1}{2}}^{n-\frac{1}{3}}\right)_{\Omega\times J} \\
& = P_{o,j_{0}}^{n-\frac{1}{3}} - P_{o,j_{0}+1}^{n}.
\end{aligned}
\label{eqn:linflu2f}
\end{equation}
Similarly, testing with $\varphi_{j_{0}+\frac{1}{2}}^{n-\frac{2}{3}}$ and $\varphi_{j_{0}+\frac{1}{2}}^{n-1}$ we get,
\begin{equation}
\begin{aligned}
\left(p_{o},\nabla \cdot \varphi_{j_{0}+\frac{1}{2}}^{n-\frac{2}{3}}\right)_{\Omega\times J} & = P_{o,j_{0}}^{n-\frac{2}{3}} - P_{o,j_{0}+1}^{n},\\
\text{and}, \\
\left(p_{o},\nabla \cdot \varphi_{j_{0}+\frac{1}{2}}^{n-1}\right)_{\Omega\times J} & = P_{o,j_{0}}^{n-1} - P_{o,j_{0}+1}^{n},
\end{aligned}
\label{eqn:linflu2c}
\end{equation}
respectively. The third term in Eqn. \eqref{eqn:watdar} involving the capillary pressure can be evaluated similarly. For the water conservation equation testing with $w_{j}^{n}$ for first term in Eqn. \eqref{eqn:watcon} we get,
\begin{equation}
\begin{aligned}
\left( \frac{\partial}{\partial t} \sum_{m=1}^{q}\sum_{i=1}^{r} \phi \rho_{w} S_{w,i}^{m} w_{i}^{m}, w_{j}^{n} \right)_{\Omega \times J}  & =\left( \frac{\partial}{\partial t} {(\phi \rho_{w} S_{w})}_{j}^{n} w_{j}^{n} , w_{j}^{n}\right) + \left( {(\phi \rho_{w} S_{w})}_{j}^{n}-{(\phi \rho_{w} S_{w})}_{j}^{n-1}, w_{j}^{n-1}\right)\\
& = \left( {(\phi \rho_{w} S_{w})}_{j}^{n} - {(\phi \rho_{w} S_{w})}_{j}^{n-1}\right) |E_{j}^{n-1}|.
\end{aligned}
\end{equation}
For the coarse domain $(j_{0}+\frac{1}{2})^{+}$, the above integral remains unchanged. However, for the fine domain $(j_{0}+\frac{1}{2})^{-}$ we have,
\begin{equation}
\begin{aligned}
\left(\frac{\partial}{\partial t} (\phi \rho_{w} s_{w}), w_{j_{0}}^{n-\frac{2}{3}}\right)_{(\Omega \times J)} = & \left( {(\phi \rho_{w} S_{w})}_{j_{0}}^{n-\frac{2}{3}} - {(\phi \rho_{w} S_{w})}_{j_{0}}^{n-1}\right) |E_{j_{0}}^{n-1}|, \\
\left(\frac{\partial}{\partial t} (\phi \rho_{w} s_{w}), w_{j_{0}}^{n-\frac{1}{3}}\right)_{(\Omega \times J)} = & \left( {(\phi \rho_{w} S_{w})}_{j_{0}}^{n-\frac{1}{3}} - {(\phi \rho_{w} S_{w})}_{j_{0}}^{n-\frac{2}{3}}\right) |E_{j_{0}}^{n-\frac{2}{3}}|, \\
\left(\frac{\partial}{\partial t} (\phi \rho_{w} s_{w}) , w_{j_{0}}^{n}\right)_{(\Omega \times J)} = & \left( {(\phi \rho_{w} S_{w})}_{j_{0}}^{n} - {(\phi \rho_{w} S_{w})}_{j_{0}}^{n-\frac{1}{3}}\right) |E_{j_{0}}^{n-\frac{1}{3}}|.
\end{aligned}
\end{equation}
The second term in Eqn. \eqref{eqn:watdar} is, 
\begin{equation}
\begin{aligned}
\left(\nabla \cdot \bs{u}_{w}, w_{j}^{n}\right)_{\Omega \times J} & = \left(\nabla \cdot \bs{u}_{w}, w_{j}^{n}\right)_{E_{j}^{n}} \\
& = U_{w,j+\frac{1}{2}}^{n} - U_{w,j-\frac{1}{2}}^{n}.
\end{aligned}
\end{equation}
For a fine domain element, with an edge at the the non-matching space time interface, we can write this term as,
\begin{equation}
\begin{aligned}
\left(\nabla \cdot \bs{u}_{w}, w_{j_{0}}^{n-\frac{2}{3}}\right)_{\Omega \times J} & = U_{w,j_{0}+\frac{1}{2}}^{n-\frac{2}{3}} - U_{w,j_{0}-\frac{1}{2}}^{n-\frac{2}{3}}\\
\left(\nabla \cdot \bs{u}_{w}, w_{j_{0}}^{n-\frac{1}{3}}\right)_{\Omega \times J} & = U_{w,j_{0}+\frac{1}{2}}^{n-\frac{1}{3}} - U_{w,j_{0}-\frac{1}{2}}^{n-\frac{1}{3}}\\
\left(\nabla \cdot \bs{u}_{w}, w_{j_{0}}^{n}\right)_{\Omega \times J} & = U_{w,j_{0}+\frac{1}{2}}^{n} - U_{w,j_{0}-\frac{1}{2}}^{n}
\end{aligned}
\label{eqn:lincon2f}
\end{equation}
Similarly, for a coarse domain element,
\begin{equation}
\begin{aligned}
\left(\nabla \cdot \bs{u}_{w}, w_{j_{0}+1}^{n}\right)_{\Omega \times J} & = U_{w,j_{0}+\frac{3}{2}}^{n} - U_{w,j_{0}+\frac{1}{2}}^{n-\frac{2}{3}} - U_{w,j_{0}+\frac{1}{2}}^{n-\frac{1}{3}} - U_{w,j_{0}+\frac{1}{2}}^{n}\\
\end{aligned}
\label{eqn:lincon2c}
\end{equation}
The third term in Eqn. \eqref{eqn:watdar} is easy to evaluate for both the coarse and fine domains as,
\begin{equation}
\begin{aligned}
\left(q_{w}, w_{j}^{n}\right) _{\Omega \times J} = q_{w,j}^{n} \left|E_{j}^{n}\right|.
\end{aligned}
\label{eqn:linconsource}
\end{equation}
The integral terms in the total mass conservation Eqn. \eqref{eqn:totcon} can be expanded by simply following the above treatment for the water conservation equation. We therefore curtail its description here to avoid redundancy. The first and second terms in Eqn. \eqref{eqn:expand} are approximated as,
\begin{equation}
\begin{aligned}
(\bs{u}_{\alpha}, \bs{v})  \approx (\bs{u}_{\alpha}, \bs{v})_{Q}
	 & = (\bs{u}_{\alpha}, \phi^{n}_{j+\frac{1}{2}})_{Q} \\[1ex] 
	& = \sum_{m=1}^{q}\sum_{i = 1}^{r+1} U^{m}_{\alpha,i+1/2} (\phi^{m}_{i+\frac{1}{2}}, \phi^{n}_{j+\frac{1}{2}})_{Q}
	   = \frac{h_j{} + h_{j+1}}{2\, \rvert e^{n}_{j+\frac{1}{2}} \rvert} \,U^{n}_{\alpha,j+\frac{1}{2}}, \text{and}
	    \label{eqn:aux_left_bf}		
\end{aligned}
\end{equation}
\begin{equation} \label{eqn:aux_right_bf}
(\lambda_{\alpha}\tilde{\bs{u}}_{\alpha}, \bs{v}) \approx (\lambda^{*}_{\alpha}\tilde{\bs{u}}_{\alpha},\bs{v})_{Q} 
 = \frac{h_{j}+ h_{j+1}}{2\, \rvert e^{n}_{j+\frac{1}{2}} \rvert} \, \lambda^{*,n}_{\alpha,j+\frac{1}{2}} \tilde{U}^{n}_{\alpha,j+\frac{1}{2}}
\end{equation}
, respectively. Here, is $Q$ the appropriate quadrature rule and $\lambda^{*,n}_{\alpha,j+\frac{1}{2}}$ is the upwind mobility defined as,
\begin{equation}\label{eqn:upw_mob}
\lambda^{*,n}_{\alpha,j+\frac{1}{2}} = \rho^{*,n}_{\alpha,j+\frac{1}{2}} \dfrac{k^{r\alpha,*}_{j+\frac{1}{2}}}{\mu_{\alpha}} =
\begin{cases}
\dfrac{1}{2\mu_{\alpha}} \left(\rho^{n}_{\alpha,j} + \rho^{n}_{\alpha,j+1}\right) k_{r\alpha}(S^{n}_{\alpha,j}), \quad & \text{~if ~} \tilde{U}^{n}_{\alpha,j+\frac{1}{2}} > 0,\\[3ex]
\dfrac{1}{2\mu_{\alpha}}  \left(\rho^{n}_{\alpha,j} + \rho^{n}_{\alpha,j+1}\right) k_{r\alpha}(S^{n}_{\alpha,j+1}),  & \text{~otherwise}.
\end{cases}
\end{equation}
This gives us a non-linear, algebraic system of equations in pressure ($P_{o}$), saturation ($S_{w}$), and Darcy ($U_{\alpha}$) and auxiliary ($\tilde{U}_{\alpha}$) flux unknowns. We now have a system of non-linear algebraic equations that can be linearized using Newton's method (or its variants) to form a linear system of equations that can be solved to obtain the solution. Please note that he flux unknowns will be eliminated using multiple Schur complements of the linearized system of equations resulting in a reduced system in pressure ($P_{o}$) and saturation ($S_{w}$) unknowns only, later in the solution algorithm section.

\section{Solution Algorithm}\label{sec:algo}
\begin{figure}[H]
\begin{center}
\includegraphics[width=7cm]{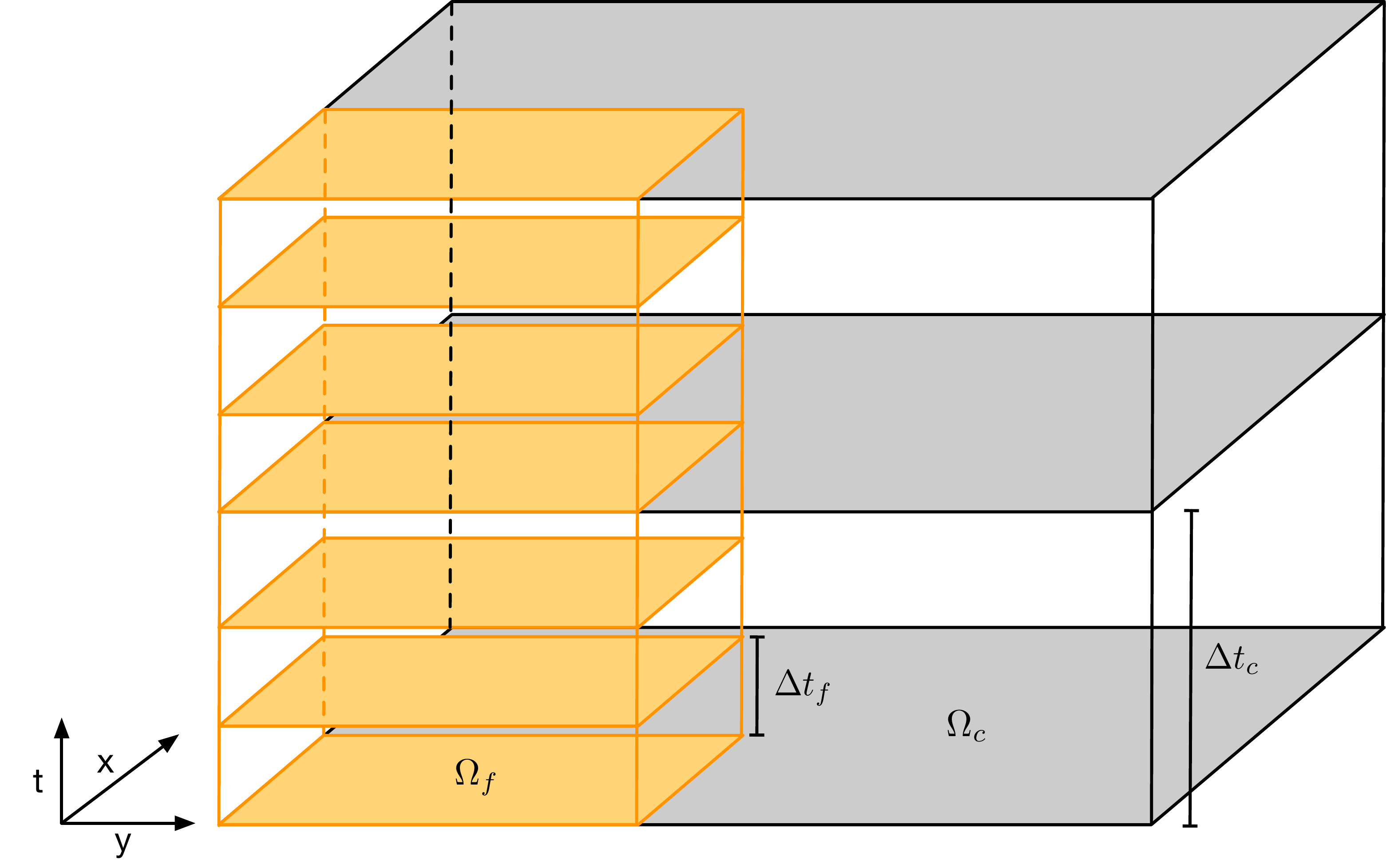}\quad
\caption{Space-time domain decomposition schematic and monolithic system construction.}
\label{fig:scheme}
\end{center}
\end{figure}

In this section, we present a space-time monolithic solver for the algebraic system obtained in the previous section. We introduce a definition for matching times as the location along the time dimension where the coarse and fine domain boundaries  overlap such that $ \frac{\Delta t_{c}}{\Delta t_{f}} = l \in \mathbb I $, where $l$ is a non-zero integer, $\Delta t_{c}$ the coarse time-step, and $\Delta t_{f}$ the fine time-step. Figure \ref{fig:scheme} shows a schematic of the space-time domain decomposition used for the construction of the monolithic system with $l = 3$. Here, the orange and grey shaded regions represent the fine and coarse, domains respectively both in space and time matching with the description of identifier 1 as discussed in Table 1 before. The schematic shown in Figure \ref{fig:scheme} and the choice of $l=3$ is kept the same in the fully discrete formulation section above to easily describe the relationship between discretization and the resulting linear system sparsity. However, as we will see later in the results section a more general setting can also be treated with more general domain decomposition using a combination of delta change in space and time as well as explicit error estimators as discussed before.
Now let us turn our attention to the conventional, time-marching solution algorithms for solving PDEs that rely upon previous time-step solution in order to determine the next time-step solution. Such an approach proves restrictive for our space-time adaptive mesh refinement framework with different time-step and grid block sizes in different subdomains of the reservoir. We therefore rely upon a time-concurrent, monolithic solution algorithm that constructs are larger system where fine and coarse domains (in time and irrespective of coarse and fine in space) are solved simultaneously over a prescribed coarse time domain. Let us denote this coarse time-step or more precisely the matching times with $\Delta T$. For the sake of simplicity and to present a coherent argument in line with the fully discrete variational formulation described before, we choose this matching time-step to be the coarse time-step size $\Delta T = \Delta t_{c}$. Again, this choice is not restrictive and a larger matching time-step size ($\Delta T>\Delta t_{c}$) can be chosen without any loss in generality. This monolithic system is then constructed to solve the algebraic equations in the space-time unknowns over the coarse and fine domains for each matching time-step ($\Delta T = \Delta t_{c}$) as follows,
\[
\left[ 
\begin{array}{c@{}c@{}c@{}c@{}c}
 & \vdots & \vdots & \vdots &  \\
   & \textcolor{red}{temporal} & spatial & &  \\
  \hdots & \left[\begin{array}{ccc}
                      \mathcal{T}^{n-\frac{1}{3}}_{pp}&0\\ 
                      0&\mathcal{T}^{n-\frac{1}{3}}_{ss}\\
                      \end{array}\right] 
                      & \left[\begin{array}{ccc}
                      \mathcal{S}^{n-\frac{2}{3}}_{pp}  &  \mathcal{S}^{n-\frac{2}{3}}_{ps} \\ 
                      \mathcal{S}^{n-\frac{2}{3}}_{sp}  & \mathcal{S}^{n-\frac{2}{3}}_{ss} \\
                      \end{array}\right]
                      &\textbf{0}
                      & \hdots \\    
& & \textcolor{red}{temporal} & spatial &   \\     
  \hdots & \textbf{0}
                      & \left[\begin{array}{ccc}
                       \mathcal{T}^{n-\frac{2}{3}}_{pp}&0\\ 
                      0&\mathcal{T}^{n-\frac{2}{3}}_{ss}\\
                      \end{array}\right] 
                      & \left[\begin{array}{ccc}
                      \mathcal{S}^{n}_{pp}  &  \mathcal{S}^{n}_{ps} \\ 
                      \mathcal{S}^{n}_{sp}  & \mathcal{S}^{n}_{ss} \\
                      \end{array}\right]
                      & \hdots \\
                       & \vdots & \vdots & \vdots &
                       \end{array}\right]
\left[
\begin{array}{c@{}}
\vdots\\
\delta x_{i-1}^{n-\frac{1}{3}}\\
\delta x_{i}^{n-\frac{1}{3}}\\
\delta x_{i+1}^{n-\frac{1}{3}}\\
\delta x_{i-1}^{n-\frac{2}{3}}\\
\delta x_{i}^{n-\frac{2}{3}}\\
\delta x_{i+1}^{n-\frac{2}{3}}\\
\delta x_{i-1}^{n}\\
\delta x_{i}^{n}\\
\delta x_{i+1}^{n}\\
\vdots
\end{array}
\right]
=
\left[
\begin{array}{c@{}}
\vdots\\
-\mathcal{R}_{i-1}^{n-\frac{1}{3}}\\
-\mathcal{R}_{i}^{n-\frac{1}{3}}\\
-\mathcal{R}_{i+1}^{n-\frac{1}{3}}\\
-\mathcal{R}_{i-1}^{n-\frac{2}{3}}\\
-\mathcal{R}_{i}^{n-\frac{2}{3}}\\
-\mathcal{R}_{i+1}^{n-\frac{2}{3}}\\
-\mathcal{R}_{i-1}^{n}\\
-\mathcal{R}_{i}^{n}\\
-\mathcal{R}_{i+1}^{n}\\
\vdots
\end{array}
\right].
\] 

Here, $\mathcal{T}$ and $\mathcal{S}$ here represent the temporal and spatial submatrices in the linearized system with subscripts $ps$ and $sp$ indicating pressure-saturation coupling submatrices. Further, $\delta x^{m}_{j} := [\delta P^{m}_{j} \text{ } \delta S_{j}^{m}]^{T}$ is the vector corresponding to the pressure and saturation unknowns and $\mathcal{R}_{j}^{m} := [R_{j,P}^{m} \text{ } R_{j,S}^{m}]$  is the residual vector corresponding to the total and water phase mass conservation equations. Note that here the flux unknowns $\delta U$ and $\delta \tilde{U}$ are eliminated by taking a Schur-complements of the original linear algebraic system in pressure, flux, and saturation unknowns. Also note that the superscripts indicating time vary as $n-1/3$, $n-2/3$, and $n-1$ in agreement with our choice of $\Delta T = \Delta t_{c}$ and the ratio between fine and coarse time-step sizes of $l = \frac{\Delta t_{c}}{\Delta t_{f}} = 3$. This monolithic construction therefore gives us a time-concurrent solution scheme where different time-step sizes in different subdomains of the reservoir can be easily resolved.

Let us now discuss the sparsity of the spatial and temporal matrices in this monolithic system that allows us to relate this scheme to the well known finite difference scheme. For the matching grid case, the spatial sub-matrix here has a known sparsity pattern of three, five, and seven non-zero diagonals for one, two, and three spatial dimensions, respectively similar to the well known finite difference scheme. For the non-matching grid case, the spatial sparsity pattern alters as for the original EVMFEM scheme in space \cite{WheelerEV}. The temporal sub-matrix is always diagonal for our choice of discretization in time. Since this discretization scheme is closely related to the backward Euler scheme \citep{arbogast15}, the solution scheme presented here is fully-implicit in the space-time unknowns. In fact, the components required to develop and use our proposed space-time adaptive framework are already available in all finite-difference discretization based numerical reservoir flow simulators.

\begin{figure}[H]
\begin{center}
\includegraphics[width=7cm]{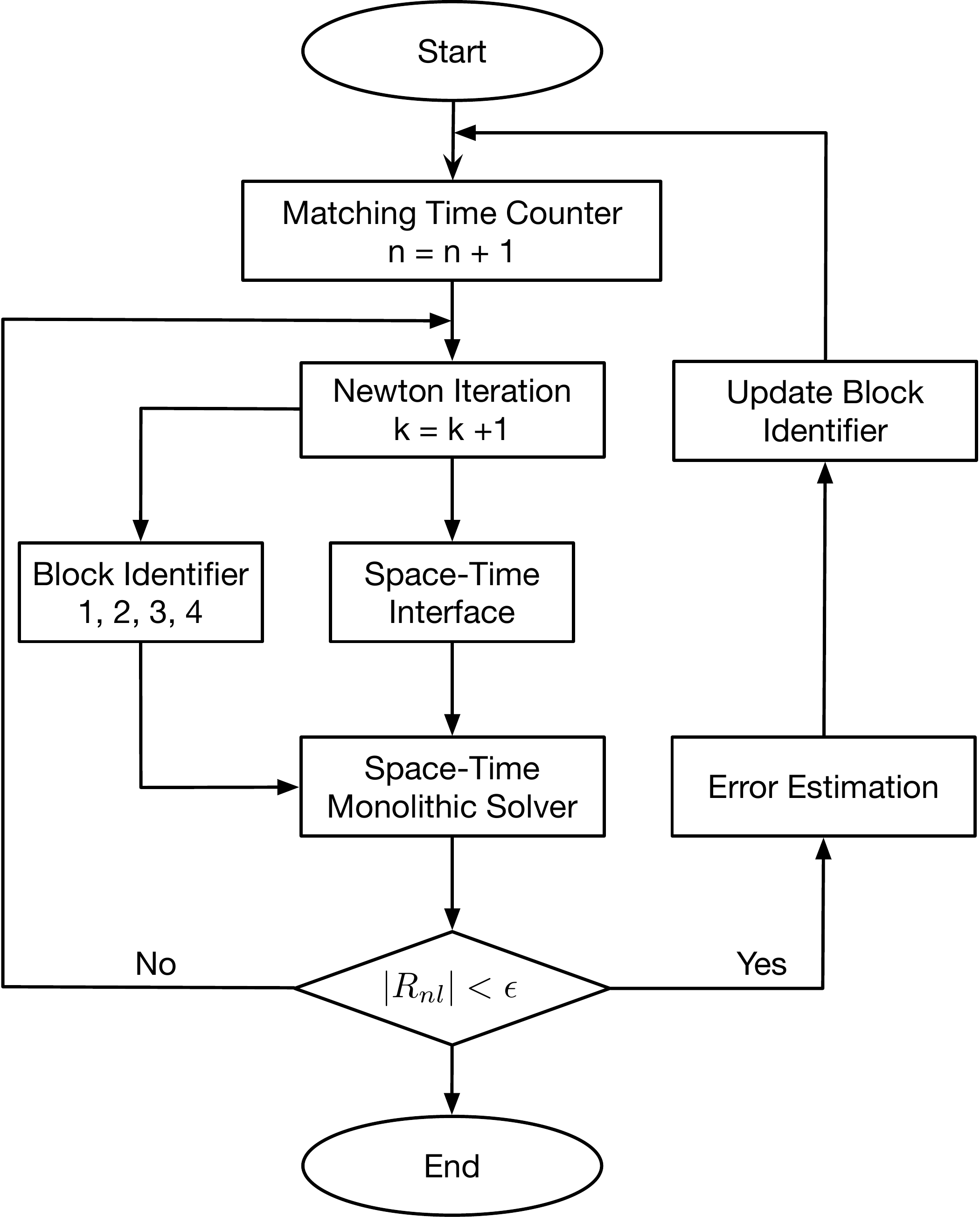}
\caption{Solution algorithm for space-time monolithic system construction and non-linear solver.}
\label{fig:flow}
\end{center}
\end{figure}

We now integrate the space-time adaptive framework using a combination of explicit error estimates (normalized non-linear residuals) and delta change in quantity of interest $Q$ to adaptively reassign identifiers (described in Table 1) during simulation runtime and decompose the space-time domain into subdomains with different spatial mesh refinements and time-step sizes to balance computational loads. Figure \ref{fig:flow} shows a flowchart for this integrated space-time adaptive framework. Here, $n$ and $k$ are the time-step and Newton iteration counters, respectively. We require that the max norm of the non-linear $\lvert R_{nl} \rvert$ residuals be less than a desired tolerance $\epsilon$ as the criteria for non-linear system convergence. We initially rely upon direct-solvers for the purpose of testing and benchmarking the solution algorithm. However, the monolithic space-time solver allows us to utilize parallel in time linear solvers and preconditioners presented in \citep{falgout14} with relative ease. This also renders us a massively parallel, time-concurrent, framework for solving general, sub-surface, non-linear flow and reactive transport problems and will be presented in a future work.

\section{Numerical results}\label{sec:results}
In this section we present two numerical experiments for the two-phase flow problem using our prototype, serial implementation of this framework. The first numerical result shows two-phase flow problem with static, space-time mesh refinements in different subdomains. Here, a fine spatial grid with small time-steps are used in the vicinity of the injection well with a coarse spatial grid and large time-steps away from the injection well. The second numerical result uses the above framework for dynamic space-time mesh refinement where spatial and temporal mesh changes as the saturation front evolves. A normalized non-linear residual is used to assign subdomain identifiers as discussed before in Table 1. 

\begin{figure}[H]
\begin{center}
\includegraphics[width=7cm,trim=0cm 0cm 0cm 0cm, clip]{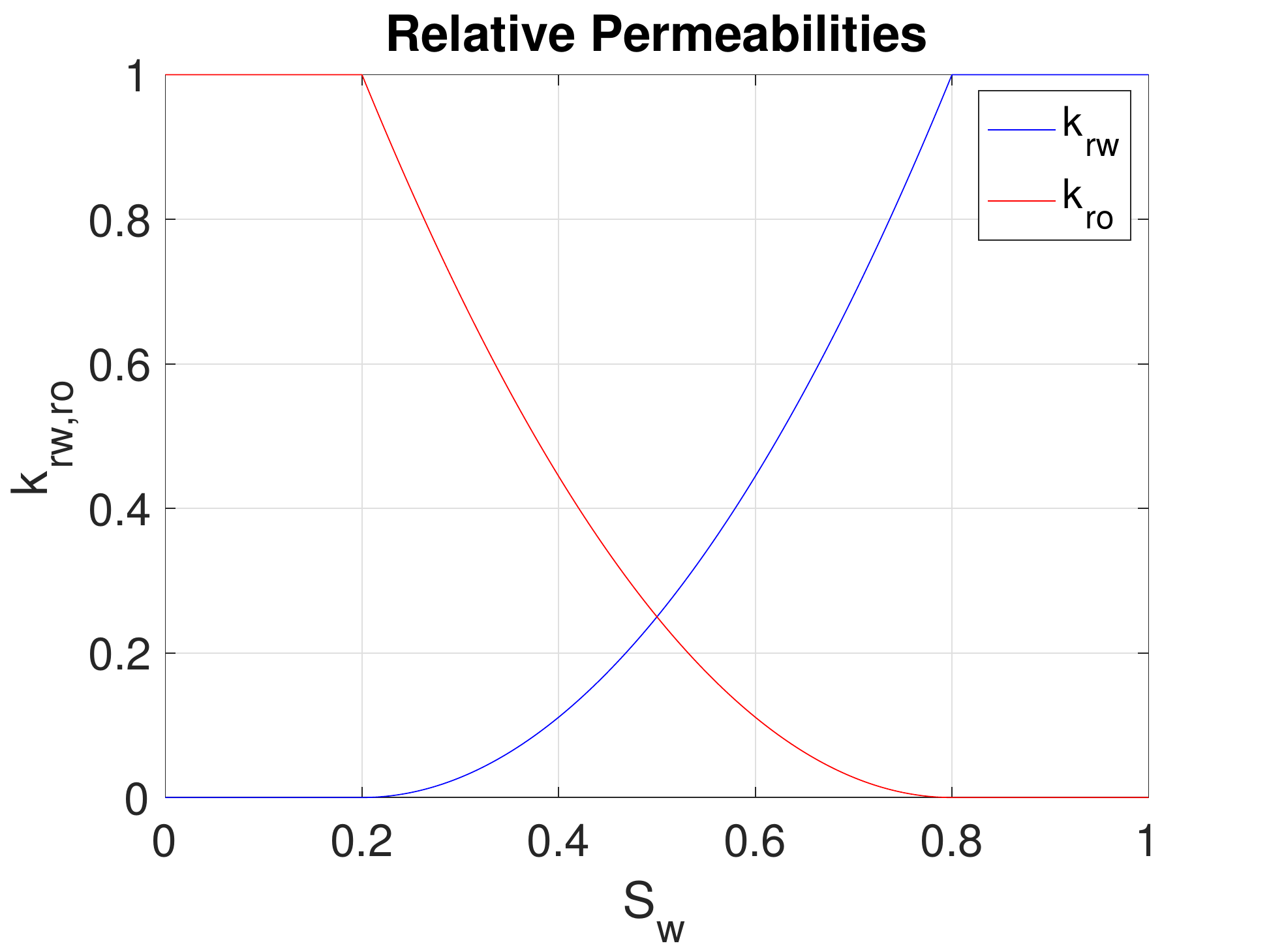}
\includegraphics[width=7cm,trim=0cm 0cm 0cm 0cm, clip]{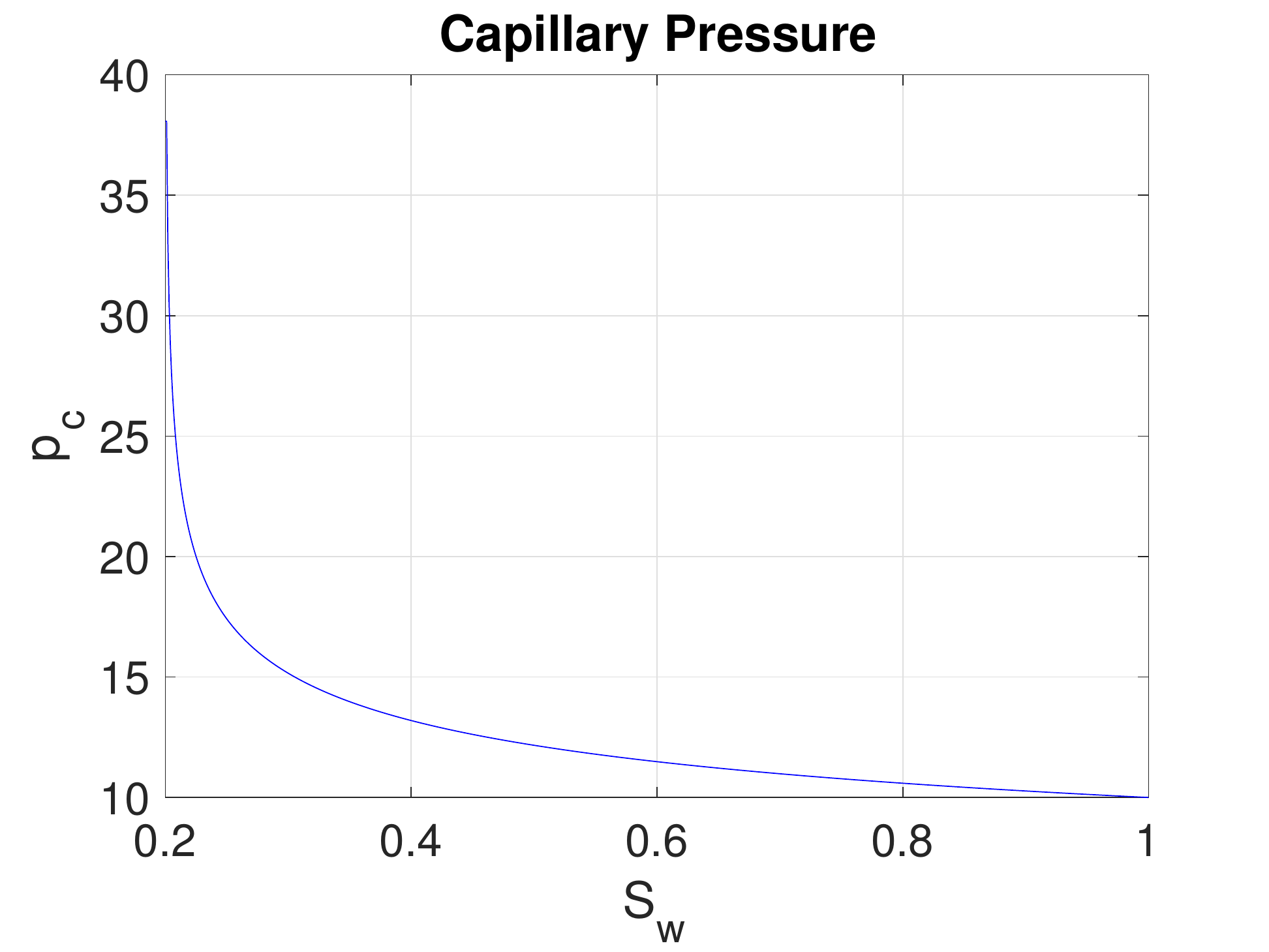}
\caption{Relative permeability (left) and capillary pressure (right) curves for the two-phase flow problem.}
\label{fig:relcap}
\end{center}
\end{figure}

For the two-phase flow model, the fluid and reservoir properties are adapted from the SPE10 \citep{spe10} dataset with an assumed homogeneous, spatial distribution for porosity of 0.2. The oil and water phase compressibilities are taken to be 1$\times$10$^{-4}$ and 3$\times$10$^{-6}~psi^{-1}$ , respectively, and densities 53$lb/ft^{3}$ and 64$lb/ft^{3}$, respectively. Further, the fluid viscosities are assumed to be 3 and 1 $cP$ for the oil and water phases, respectively. Additionally, a Brook's Corey model, Eqn. \eqref{eqn:brooks}, is considered for the two-phase relative permeabilities with endpoints S$_{or}$ = S$_{wirr}$ = 0.2 and $k^{0}_{ro}=k^{0}_{rw}=1.0$, and model exponents $n_{o}=n_{w}=2$. Figure \ref{fig:relcap} shows the relative permeability and capillary pressure curves as a function of saturation.
\begin{equation}
\begin{aligned}
k_{rw} &= k^{0}_{rw}\left(\frac{S_{w}-S_{wirr}}{(1-S_{or}-S_{wirr}}\right)^{n_{w}}\\
k_{ro} &= k^{0}_{ro}\left(\frac{S_{o}-S_{or}}{(1-S_{or}-S_{wirr}}\right)^{n_{o}}
\label{eqn:brooks}
\end{aligned}
\end{equation}
Further, the capillary pressure function is also defined using the Brooks-Corey model (Eqn. \eqref{eqn:brookscap}). The model parameters $P_{en,cow}$ and $l_{cow}$ are chosen to be 0.2 and 10 psi, respectively.
\begin{equation}
p_{c}(S_{w}) = P_{en,cow} \left(\frac{1-S_{wirr}}{S_{w}-S_{wirr}}\right)^{l_{cow}}
\label{eqn:brookscap}
\end{equation}

\subsection{Two Phase Flow Problem with Static Domain Decomposition}

\begin{figure}[H]
\begin{center}
\includegraphics[width=7.5cm,trim=2.5cm 1cm 2.5cm 1cm, clip]{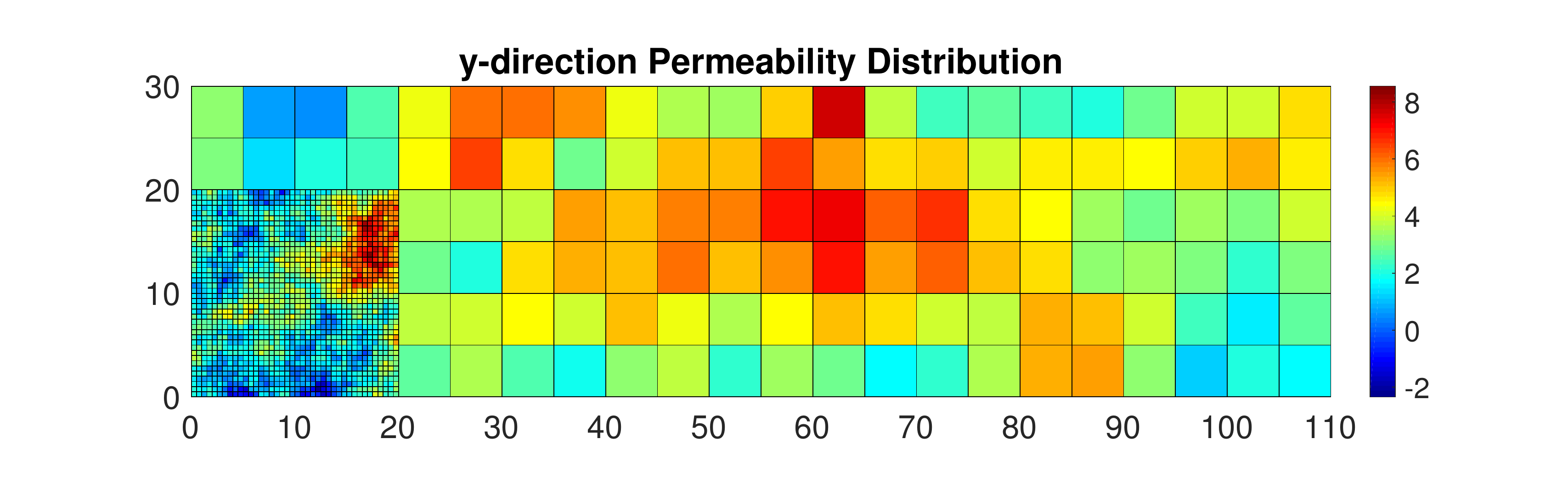}
\includegraphics[width=7.5cm,trim=2.5cm 1cm 2.5cm 1cm, clip]{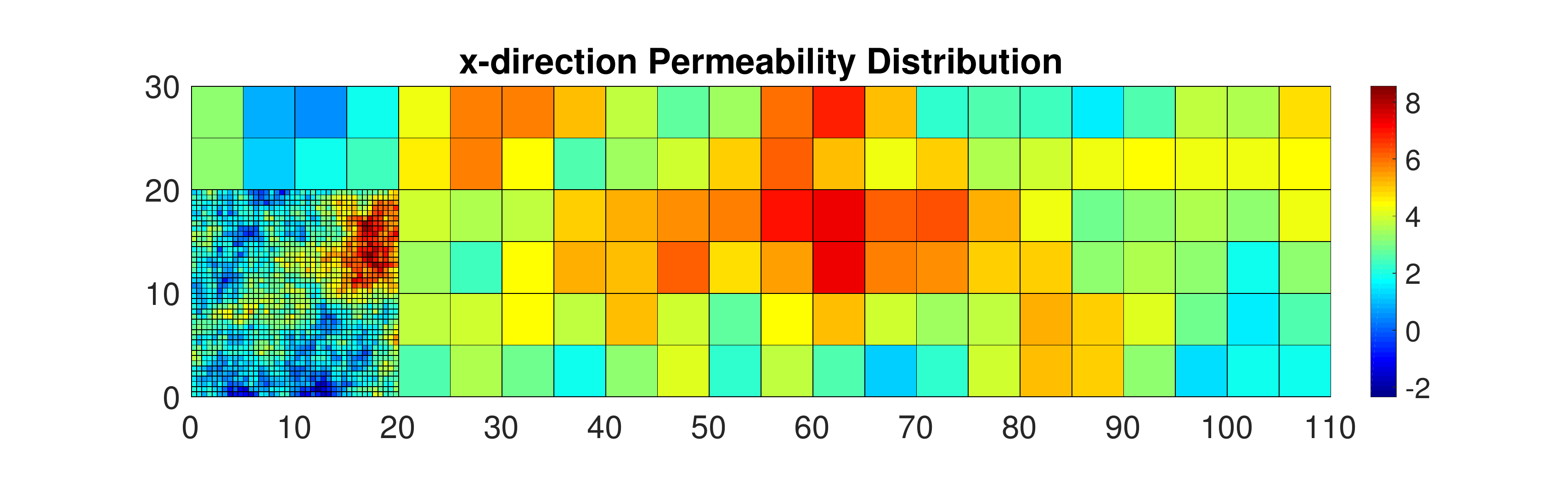}
\caption{Log-scale permeability distribution for the horizontal (left) and vertical direction (right).}
\label{fig:perml20static}
\end{center}
\end{figure}

The computational domain is 110ft $\times$ 30ft $\times$ 1ft $\times$ 40days with the fine and coarse subdomains discretized using grid elements of size 0.5ft $\times$ 0.5ft $\times$ 1ft $\times$ 1day and 5ft $\times$ 5ft $\times$ 1ft $\times$ 5days, respectively. The fine subdomain is refined by a factor of 10 and 5 times with respect to the coarse subdomain in space and time, respectively. The permeability values in the coarse subdomain are obtained by applying a local, two-scale homogenization approach given the fine scale properties from layer 20 of the 10$^{th}$ SPE comparative project dataset. The log (natural log) scale, permeability distribution in the horizontal (y) and vertical (x) directions are shown in Figure \ref{fig:perml20static}. The reader is referred to the adaptive homogenization approach in \cite{singhatce17,singhcg18} for further details on computing coarse scale effective properties.
The injection and production wells are located  at the bottom left and top right corners of the domain using rate specified injection and pressure specified production wells, respectively. The injection well is water-rate specified at 1 STB/day whereas the production well is pressure specified at 1000 $psi$. Further, the initial reservoir pressure and saturation are taken to be 1000 $psi$ and $0.2$, respectively. Figure \ref{fig:l202phsat} shows the evolution of the saturation front with time over the entire domain. The simulation was ran for a total of 40 days however we show the saturation distribution starting from 26 days up until 31 to demonstrate faster changes occurring in the fine domain compared to the coarse subdomain. The saturation distribution changes in the coarse subdomain at every coarse time-step increment (5 days) as opposed to fine subdomain where it changes at every fine time-step increment (1 day). 
\begin{figure}[H]
\begin{center}
\includegraphics[width=7.5cm,trim=2cm 0.5cm 2cm 0.5cm, clip]{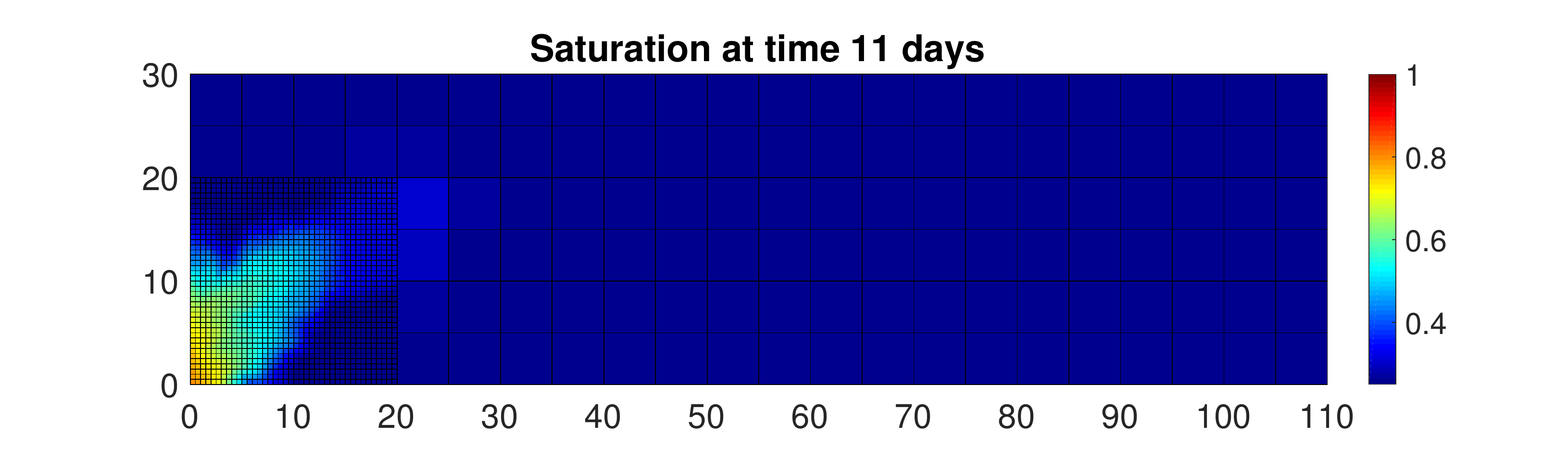}
\includegraphics[width=7.5cm,trim=2cm 0.5cm 2cm 0.5cm, clip]{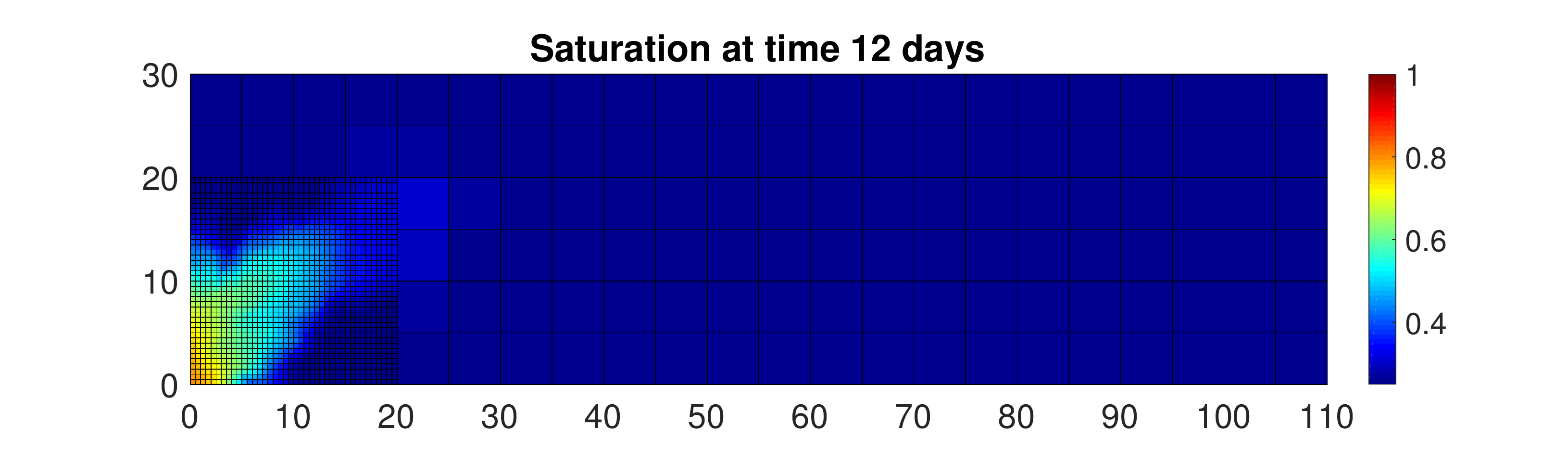}\\
\includegraphics[width=7.5cm,trim=2cm 0.5cm 2cm 0.5cm, clip]{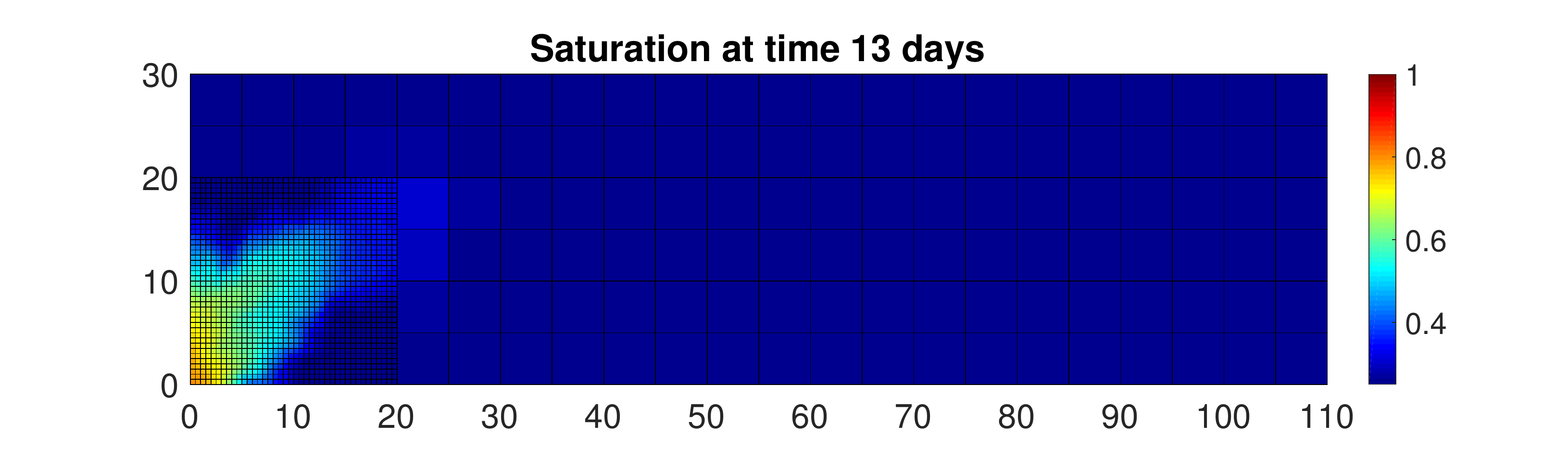}
\includegraphics[width=7.5cm,trim=2cm 0.5cm 2cm 0.5cm, clip]{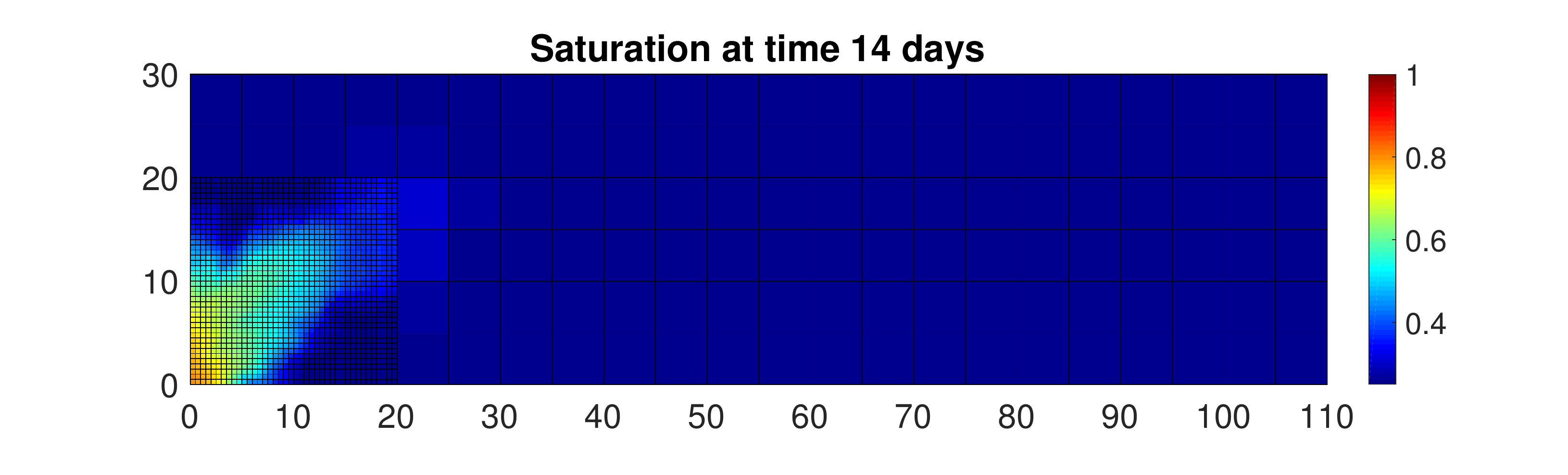}\\
\includegraphics[width=7.5cm,trim=2cm 0.5cm 2cm 0.5cm, clip]{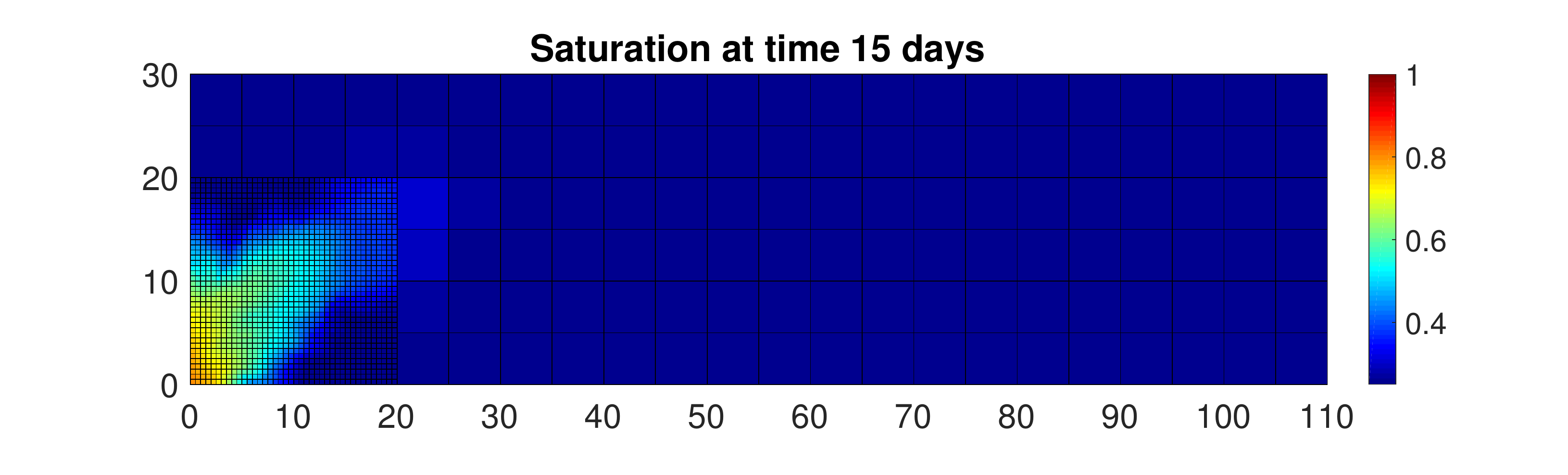}
\includegraphics[width=7.5cm,trim=2cm 0.5cm 2cm 0.5cm, clip]{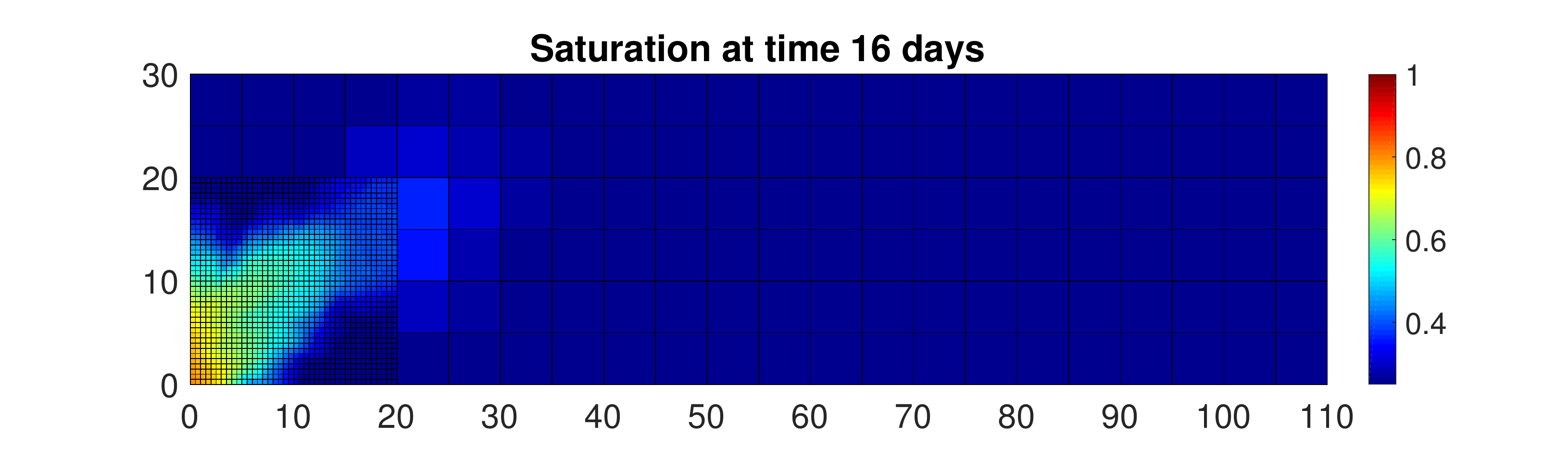}
\caption{Saturation distribution at different times.}
\label{fig:l202phsat}
\end{center}
\end{figure}
The choice of the fine subdomain, at the bottom left part of the domain, is made to capture the evolution of the saturation front starting from the injection well. It is easy to see that the non-linear functions of saturation (such as relative permeability and capillary pressure) manifest highly non-linear behavior in the region where saturation changes are large. Here, non-linear solvers such as Newton-Raphson method are marred by small time-step sizes.  An increase in time-step size often results in an increase in non-linear iterations or convergence issues.  This is further exacerbated by mesh refinement (local or global) and consequently increased computational costs. A choice of fine subdomain (space and time) in the vicinity of the saturation front allows us to not only to gain accuracy but also circumvent convergence issues associated with the non-linear solver. The computational savings are self evident since small time-step increments are only necessary in the fine subdomain as opposed to the entire domain. In fact, for a linear flow and transport problem the computational cost is directly related to the space-time degrees of freedom.

\subsection{Two Phase Flow Problem with Dynamic Domain Decomposition}
This numerical experiment demonstrates dynamic space-time mesh refinement using the framework described in the introduction section. The reservoir domain is decomposed into subdomains with identifiers listed in Table 1 based upon a combination of delta change in the quantity of interest (saturation) in space and time, and normalized non-linear residuals. The computational domain is kept the same as before (110ft $\times$ 30ft $\times$ 1ft $\times$ 100days) with the subdomains corresponding to identifier 1, 2, and 4 discretized using grid elements of size 0.5ft $\times$ 0.5ft $\times$ 0.5ft $\times$ 1day, 0.5ft $\times$ 0.5ft $\times$ 0.5ft $\times$ 4day, and 2.5ft $\times$ 2.5ft $\times$ 1ft $\times$ 4day respectively. The fine spatial and temporal mesh uses a refinement factor of 5 and 4, respectively compared coarse mesh. Please note that the refinement factors in space and time can be different  and can be altered during simulation runtime based upon our explicit estimators for better computational efficiency. However, for the sake of simplicity we use a constant refinement factor with an intention to have only two spatial and temporal scales so as to utilize the simple binary classification in Table 1. Further generalization of this approach that allows for refinement factor to be adaptively changed based upon dynamic computational loads will be presented in a future work. 

\begin{figure}[H]
\begin{center}
\includegraphics[width=7.5cm,trim=2.5cm 4.5cm 2.5cm 4.2cm, clip]{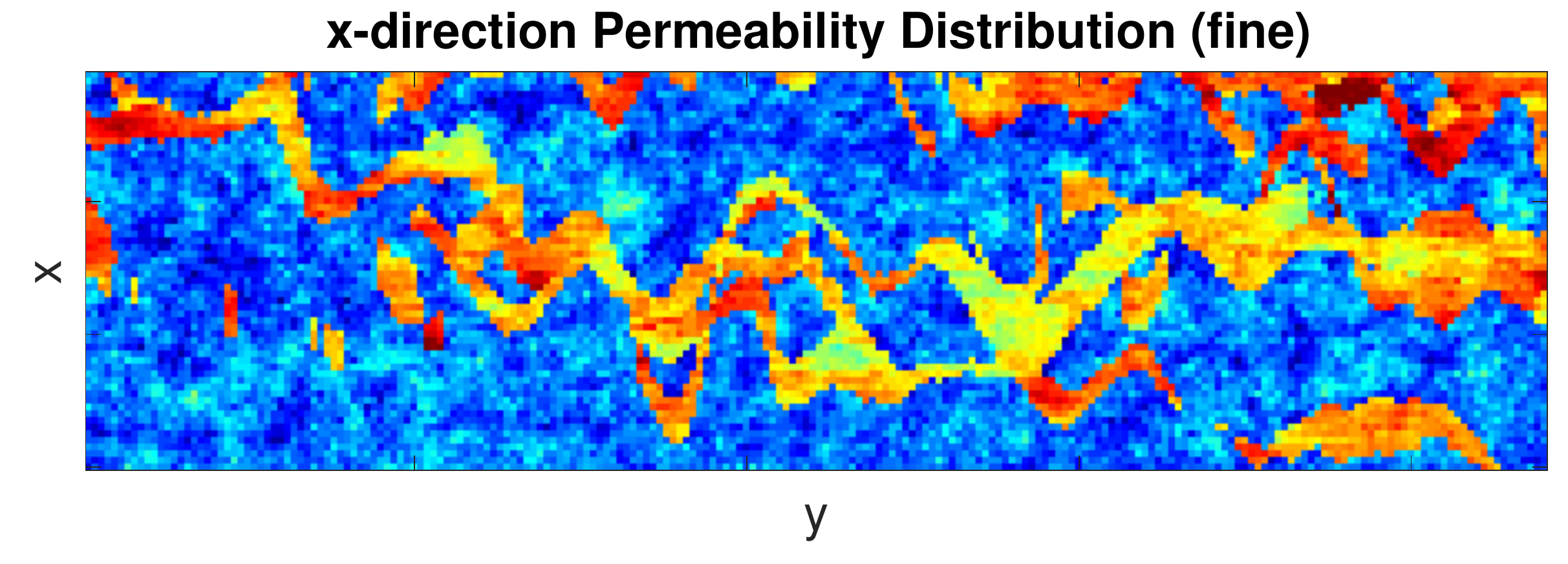}
\includegraphics[width=7.5cm,trim=2.5cm 4.5cm 2.5cm 4.2cm, clip]{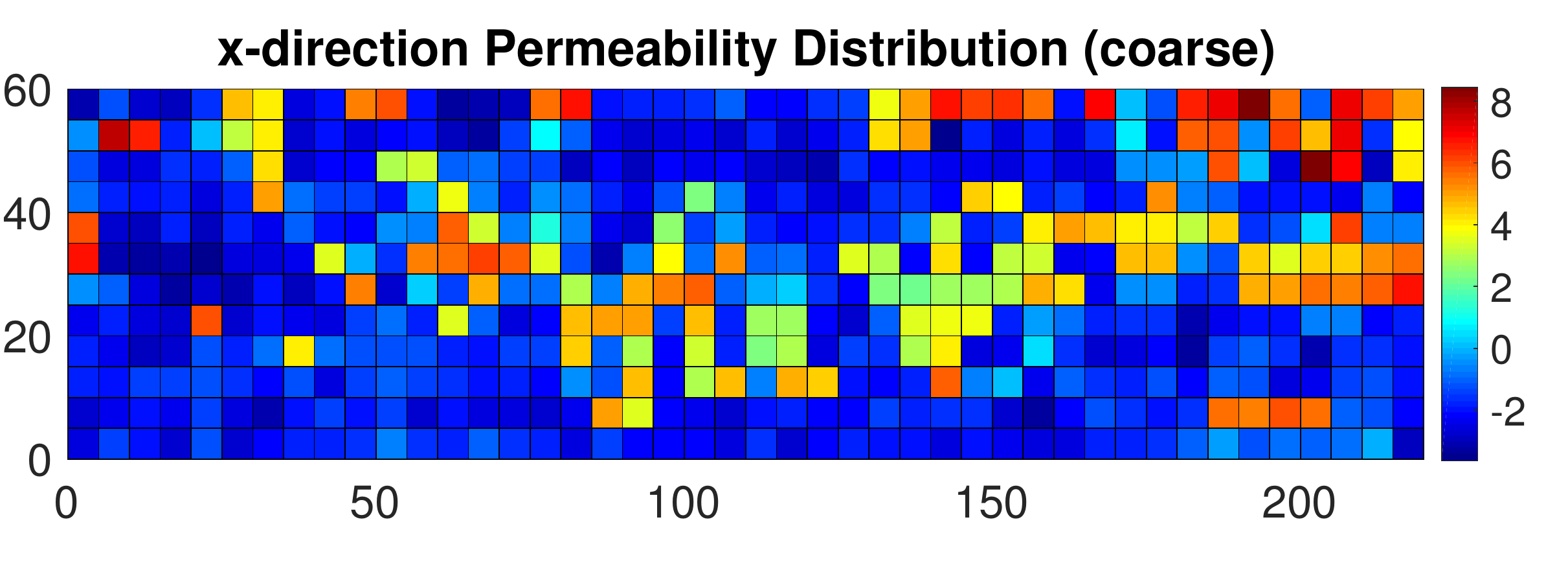}\\
\includegraphics[width=7.5cm,trim=2.5cm 4.5cm 2.5cm 4.2cm, clip]{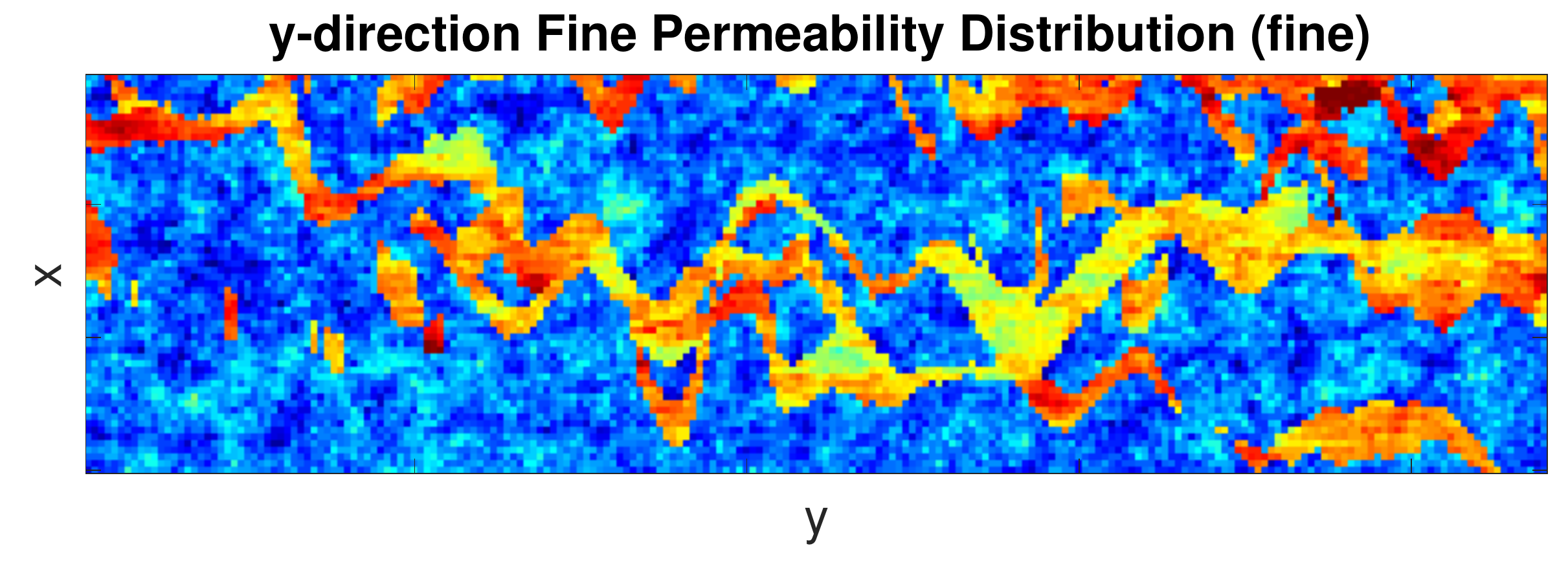}
\includegraphics[width=7.5cm,trim=2.5cm 4.5cm 2.5cm 4.2cm, clip]{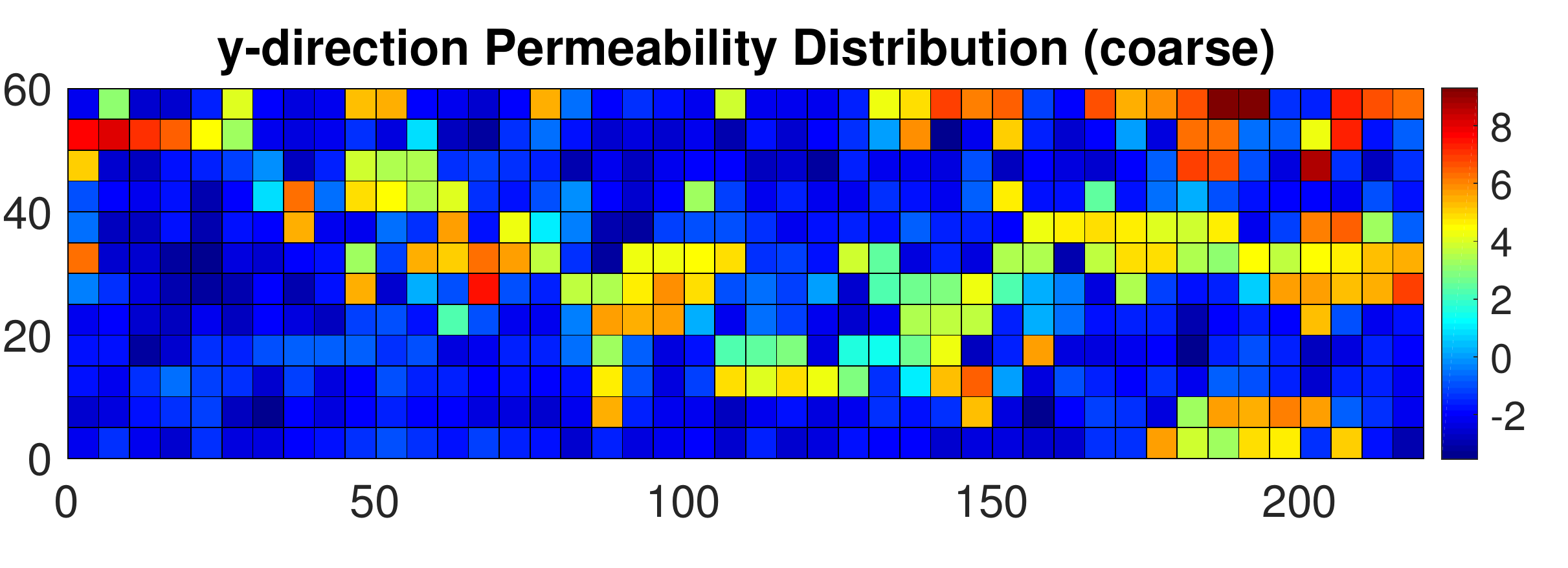}
\caption{Log-scale permeability distribution (layer 37) for the vertical (top) and horizontal directions (bottom) at the fine (left) and coarse (right) spatial scales.}
\label{fig:perml37}
\end{center}
\end{figure}

All the model parameters are kept the same as before with the exception of permeability distribution that is taken from layer 37 (instead of layer 20) of the 10$^{th}$ SPE comparative project dataset. Further, the SPE10 permeability distribution is assumed to be at a fine spatial scale (0.5ft $\times$ 0.5ft $\times$ 0.5ft) and coarse spatial scale (2.5ft $\times$ 2.5ft $\times$ 1ft) properties are obtained using a local, two-scale, numerical homogenization approach described in \cite{singhatce17,singhcg18}. Figure \ref{fig:perml37} shows the fine scale permeability distributions from the SPE10 layer 37 with coarse scale distributions from the aforementioned numerical homogenization approach. As in the previous numerical experiment injection and production wells are considered  at the bottom left and top right corners of the domain using rate specified injection and pressure specified production wells, respectively. The injection well is water-rate specified at 1 STB/day whereas the production well is pressure specified at 1000 $psi$. Further, the initial reservoir pressure and saturation are taken to be 1000 $psi$ and $0.2$, respectively.

\begin{figure}[H]
\begin{center}
\includegraphics[width=7.5cm,trim=2.5cm 0.5cm 2.5cm 0.0cm, clip]{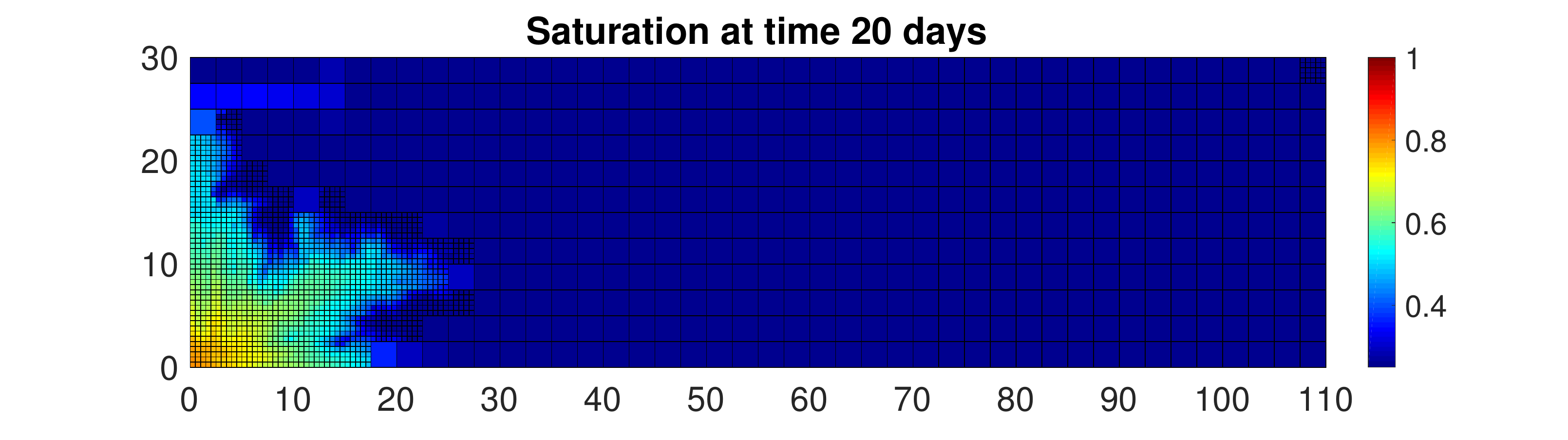}
\includegraphics[width=7.5cm,trim=2.5cm 0.5cm 2.5cm 0.0cm, clip]{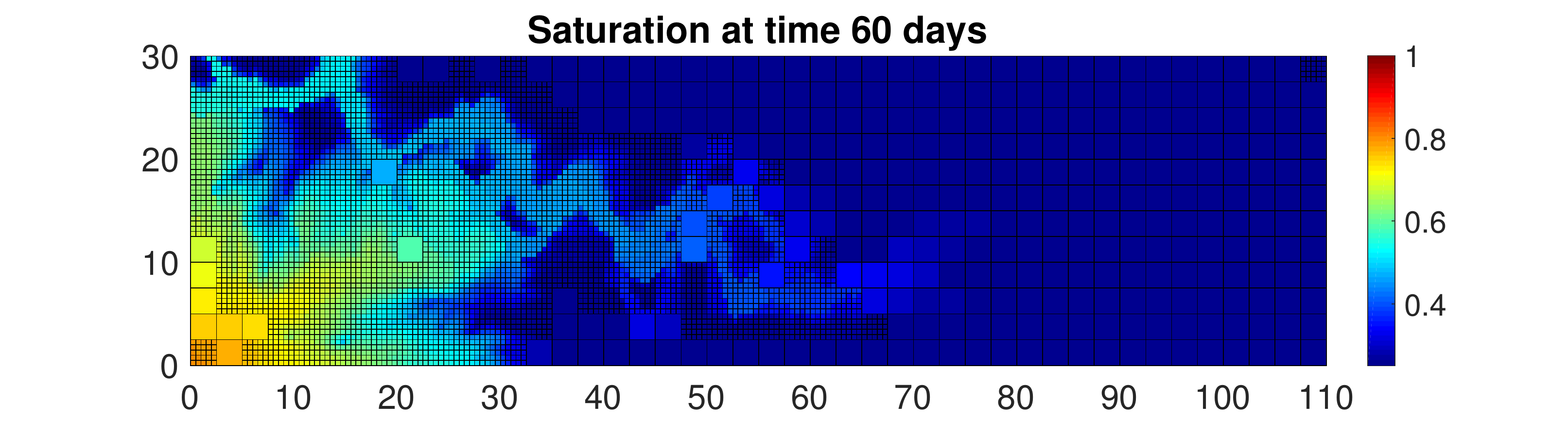}\\
 \includegraphics[width=7.4cm,trim=2.5cm 1cm 2.5cm 0.5cm, clip]{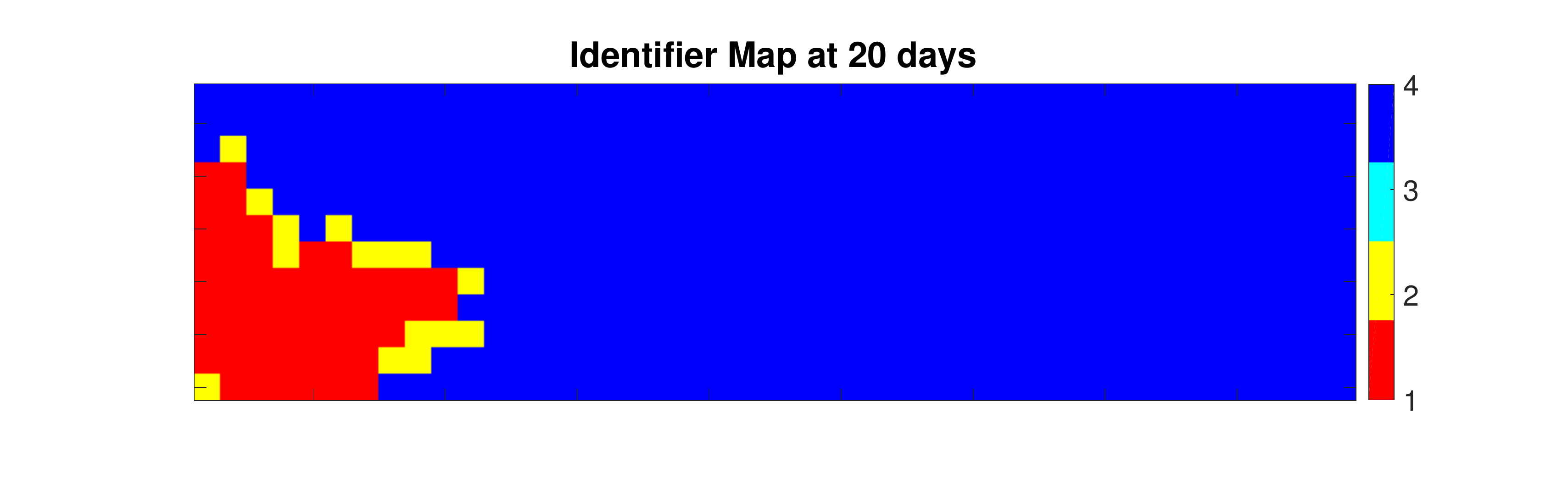} \ 
 \includegraphics[width=7.4cm,trim=2.5cm 1cm 2.5cm 0.5cm, clip]{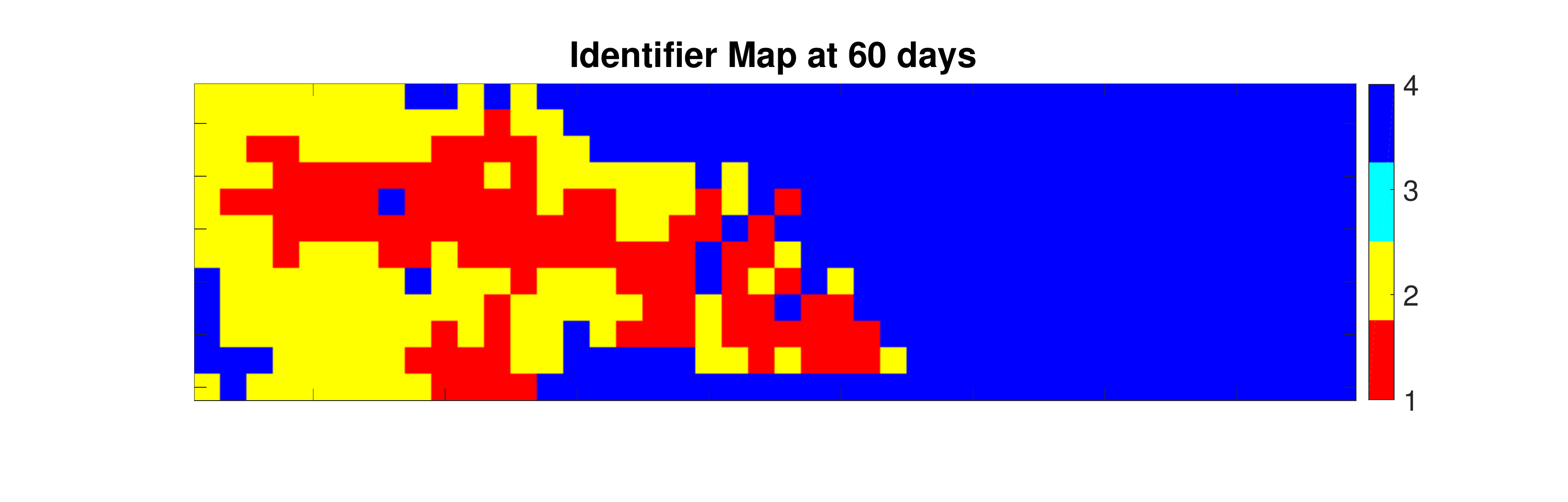} \ \ \\
\includegraphics[width=7.5cm,trim=2.5cm 0.5cm 2.5cm 0.0cm, clip]{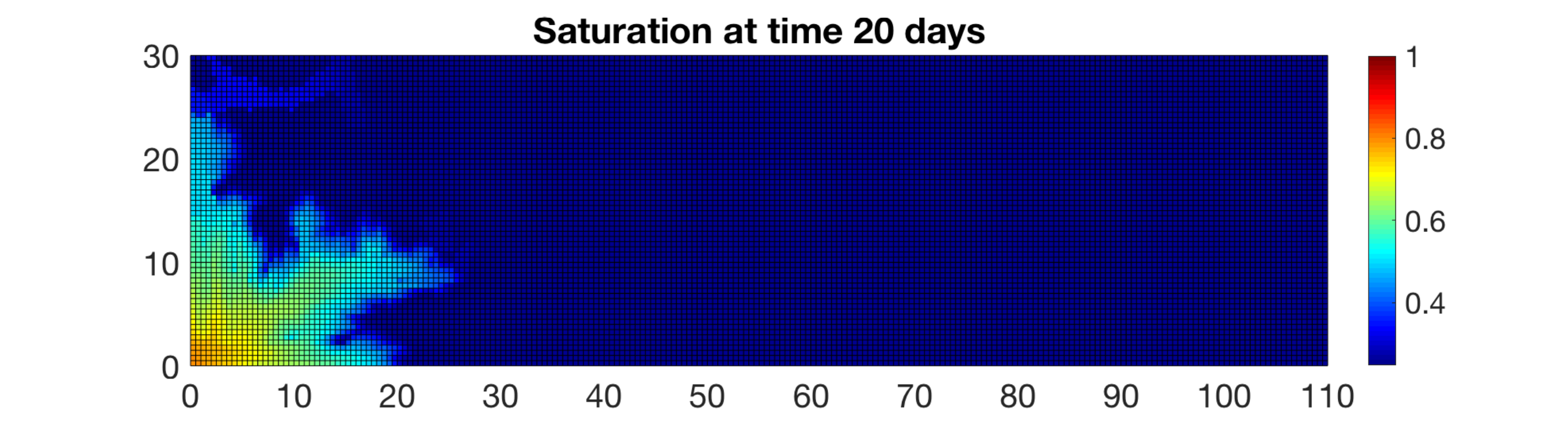}
\includegraphics[width=7.5cm,trim=2.5cm 0.5cm 2.5cm 0.0cm, clip]{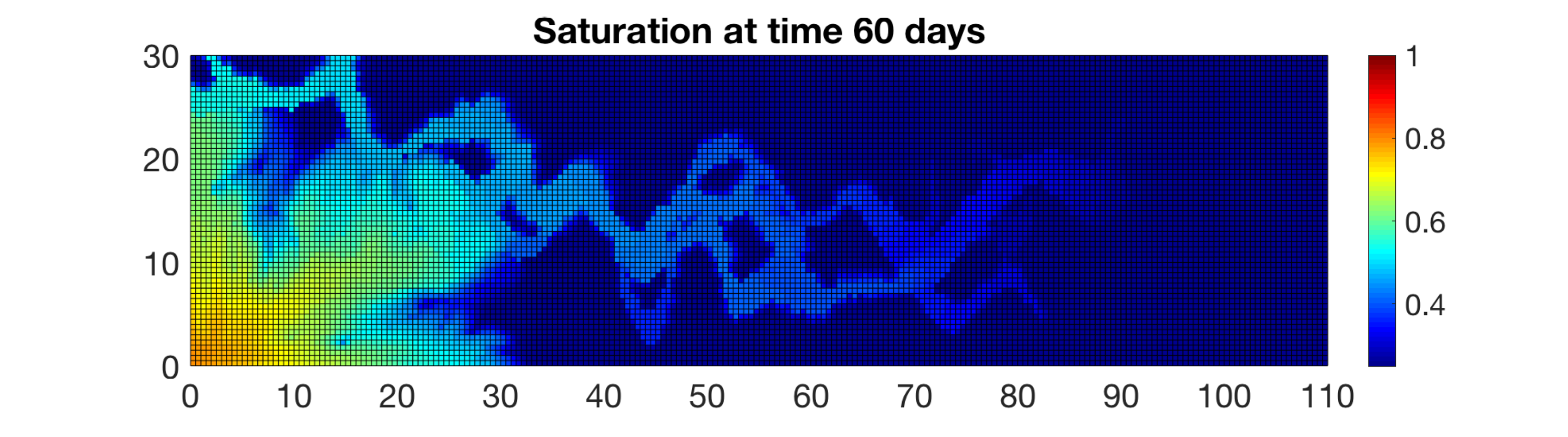}
\caption{Saturation distribution (top) and domain decomposition map (middle) for the adaptive space-time approach, and fine scale simulation comparison at 20 and 60 days (layer 37).}
\label{fig:l372phsat}
\end{center}
\end{figure}

\begin{figure}[H]
\begin{center}
\includegraphics[width=7.5cm,trim=2.5cm 0.5cm 2.5cm 0cm, clip]{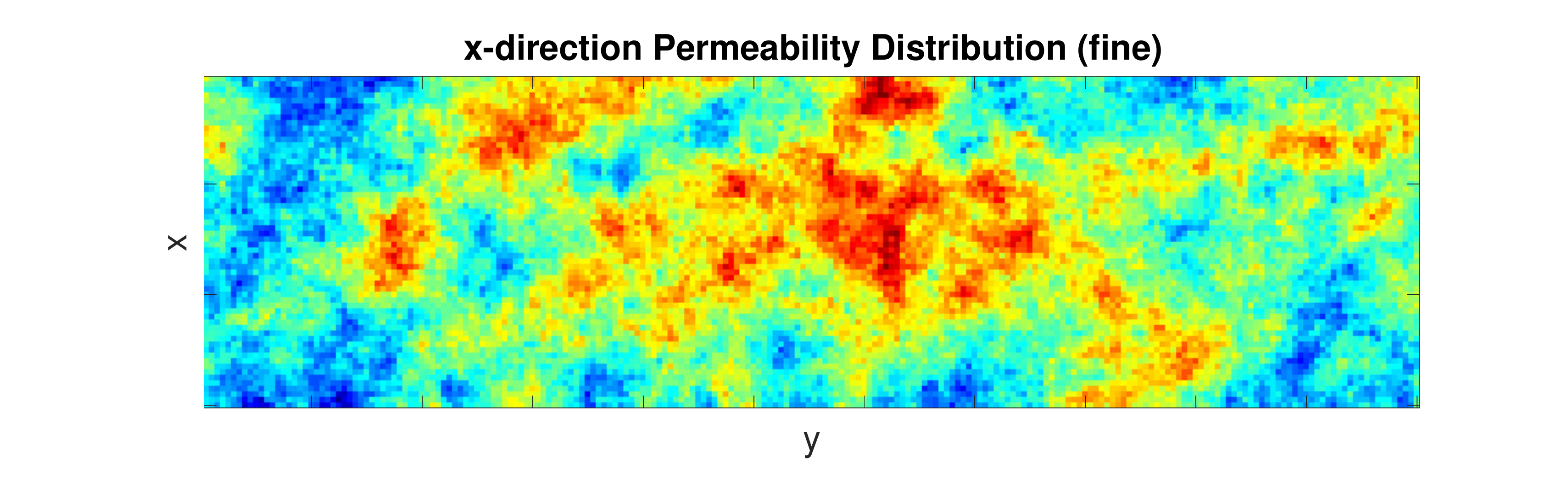}
\includegraphics[width=7.5cm,trim=2.5cm 0.5cm 2.5cm 0cm, clip]{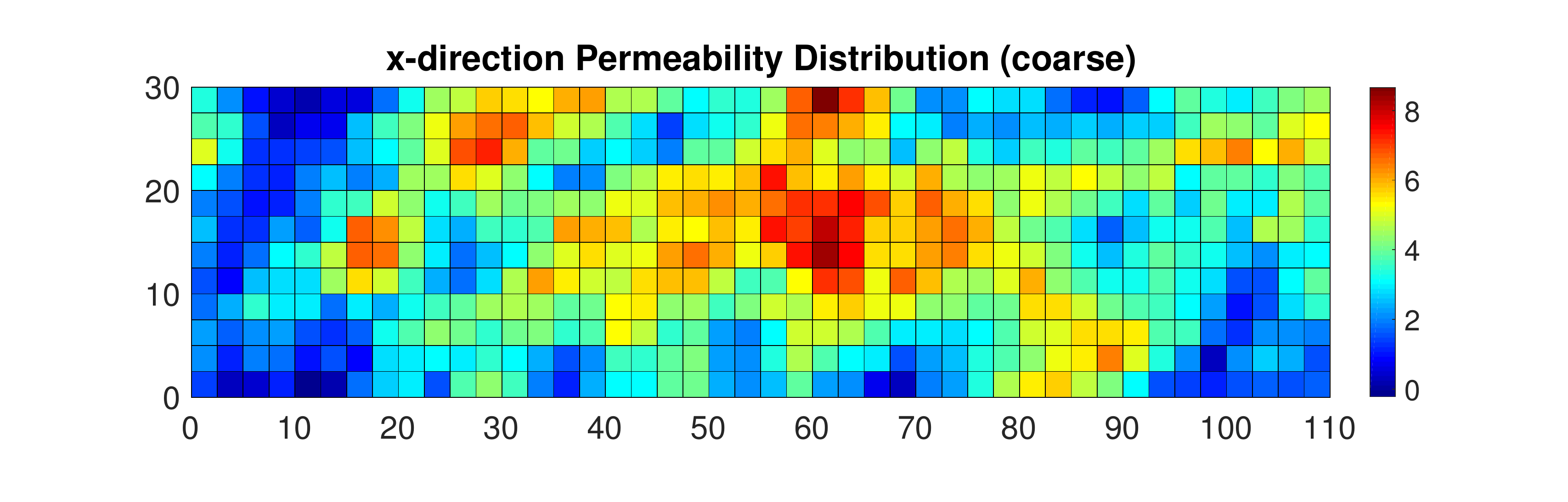}\\
\includegraphics[width=7.5cm,trim=2.5cm 0.5cm 2.5cm 0cm, clip]{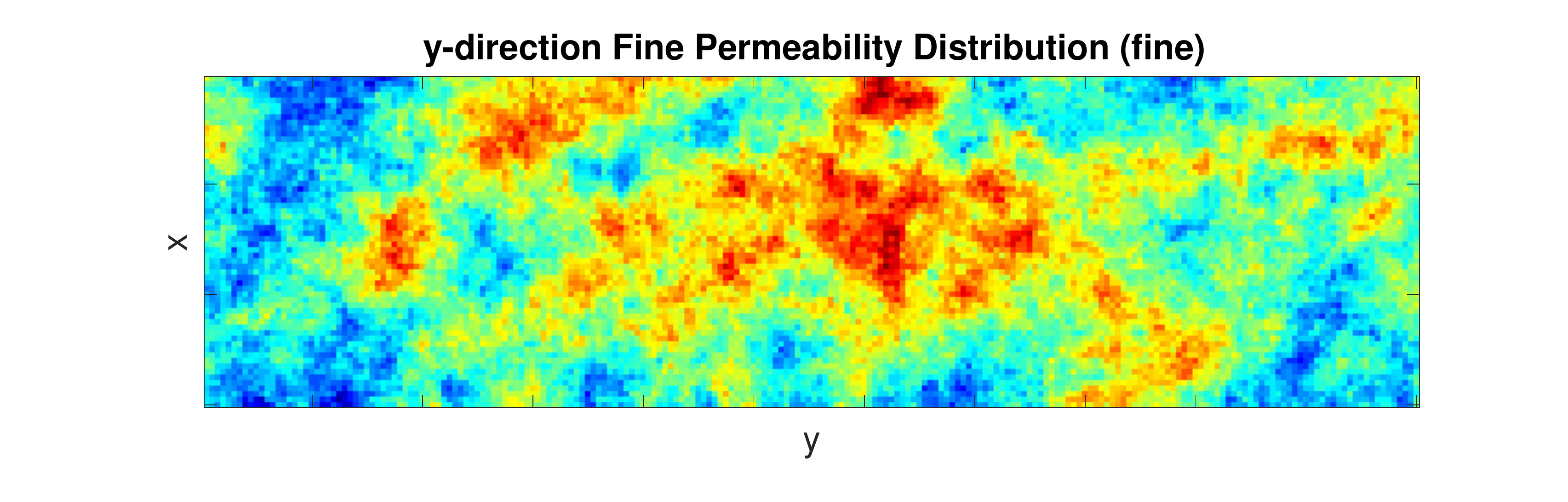}
\includegraphics[width=7.5cm,trim=2.5cm 0.5cm 2.5cm 0cm, clip]{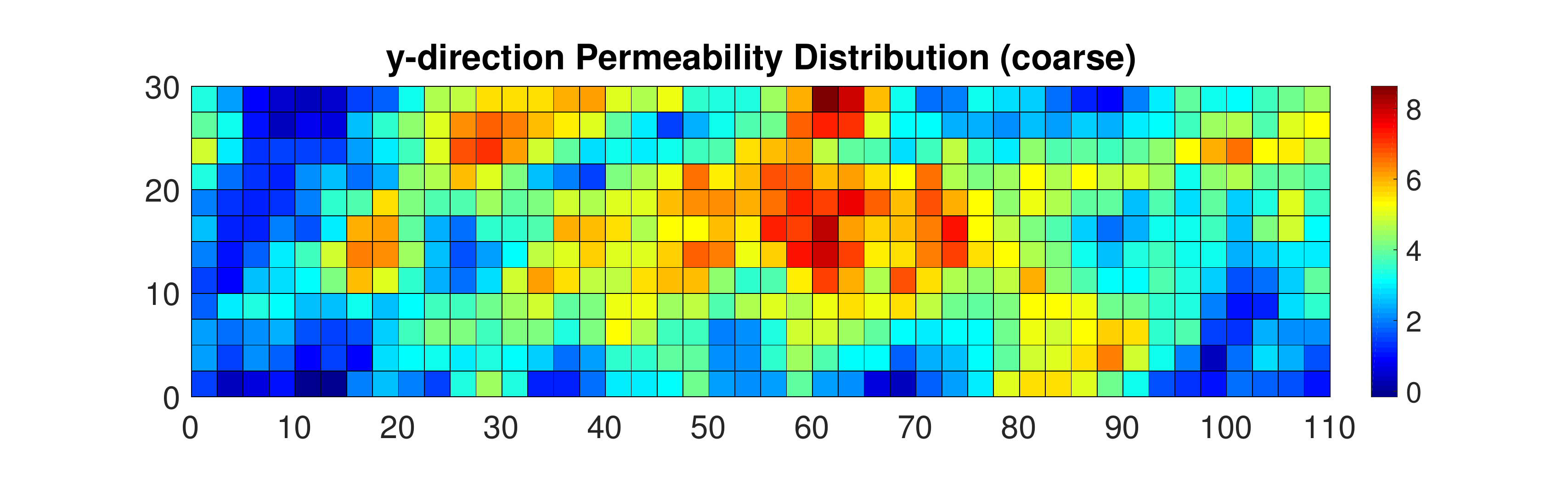}
\caption{Log-scale permeability distribution (layer 20) for the vertical (top) and horizontal directions (bottom) at the fine (left) and coarse (right) spatial scales.}
\label{fig:perml20}
\end{center}
\end{figure}

Figure \ref{fig:l372phsat} (left) shows the evolution of the saturation front with time for the adaptive space-time approach (top), the domain decomposition map (middle), and the fine scale simulation results at 20 and 60 days. The identifier map shows the classification of the reservoir domain into identifiers 1 (red), 2 (yellow), and 4 (blue) (see Table 1) using explicit estimators and delta change in the quantity of interest i.e., saturation. Please note that 1, 2, and 4 are listed here in the decreasing order of computational cost i.e., the higher the identifier the lower the computational resource required. This map is used to define spatio-temporal mesh refinement for the next coarse step $\Delta T = \Delta t_{c}$, as discussed in the solution algorithm section. The simulation was ran for a total of 100 days with dynamic spatio-temporal mesh refinements by decomposing the space-time domain into subdomains identified by classifier described in Table 1. The saturation distribution changes in the coarse subdomain at every coarse time-step increment (4 days) as opposed to fine subdomain where it changes at every fine time-step increment (1 day). 
Next, we present a similar numerical experiment with a Gaussian-like permeability distribution. We intentionally use different layers with two spatial dimensions to draw out the computational performance attributes corresponding to Gaussian and channelized permeability distributions. Figure \ref{fig:perml20} shows the fine scale permeability distributions from the SPE10 layer 20 with coarse scale distributions from numerical homogenization approach. As before, Figure \ref{fig:l202phsat} (left) shows the evolution of the saturation front with time for the adaptive space-time approach (top), the domain decomposition map (middle), and the fine scale simulation results at 20 and 60 days. 

\begin{figure}[H]
\begin{center}
\includegraphics[width=7.5cm,trim=2.5cm 1cm 2.5cm 0.5cm, clip]{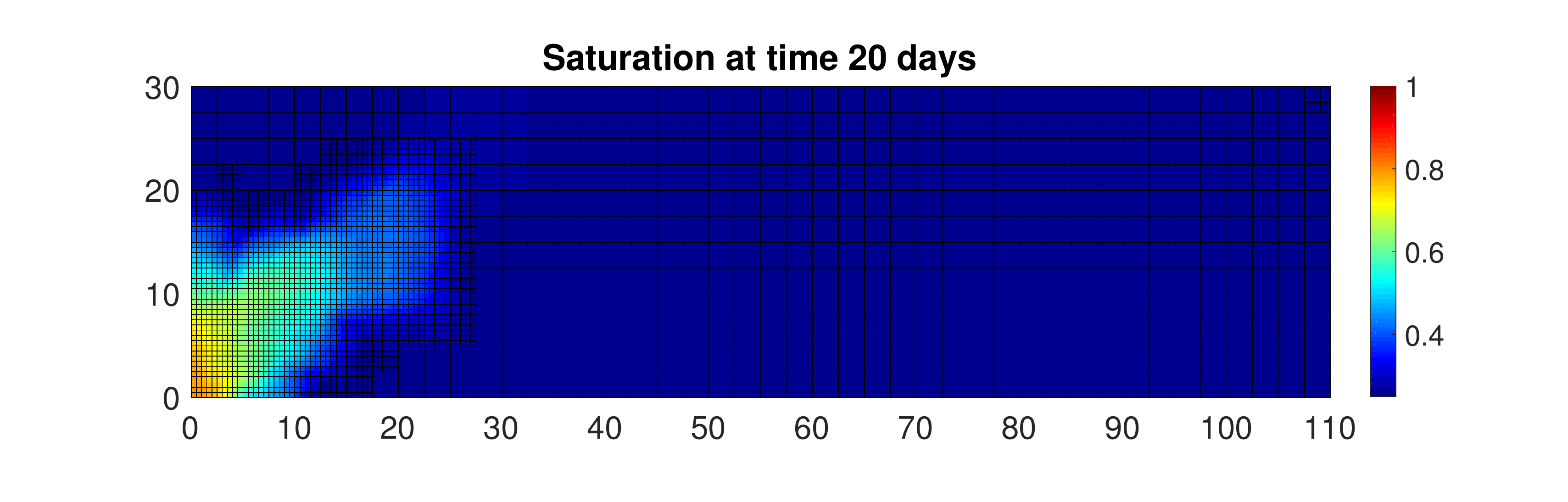}
\includegraphics[width=7.5cm,trim=2.5cm 1cm 2.5cm 0.5cm, clip]{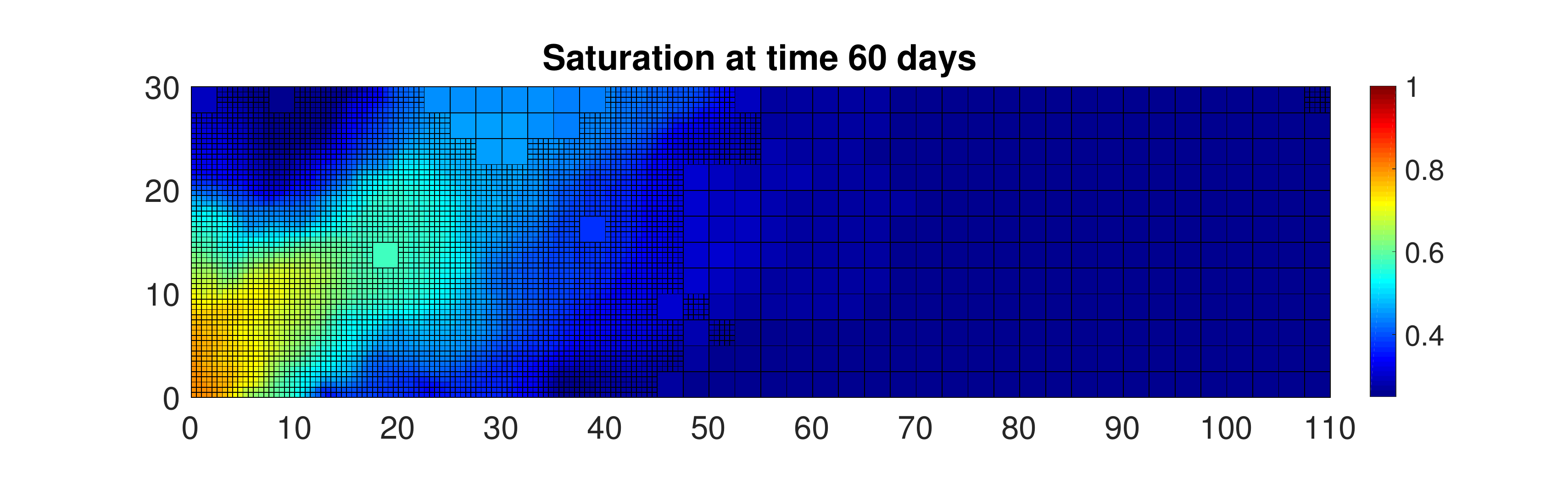}\\
 \includegraphics[width=7.3cm,trim=2.5cm 2cm 2.5cm 0.5cm, clip]{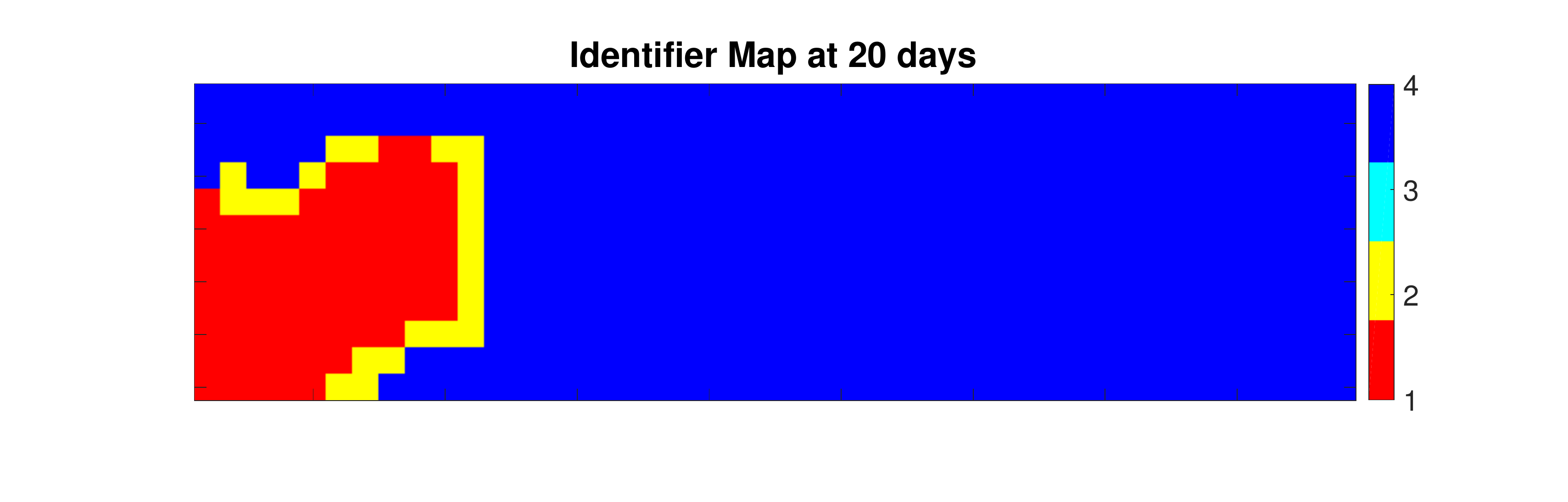} \ \
 \includegraphics[width=7.3cm,trim=2.5cm 2cm 2.5cm 0.5cm, clip]{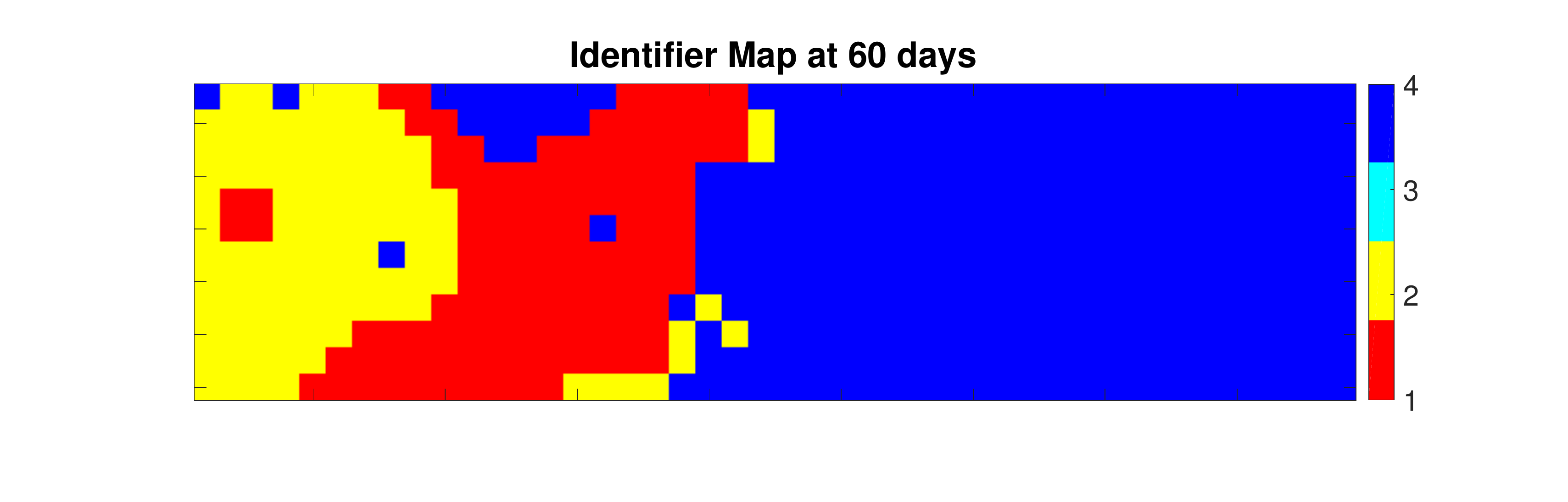} \ \ \ \\
\includegraphics[width=7.5cm,trim=2.5cm 1cm 2.5cm 0.5cm, clip]{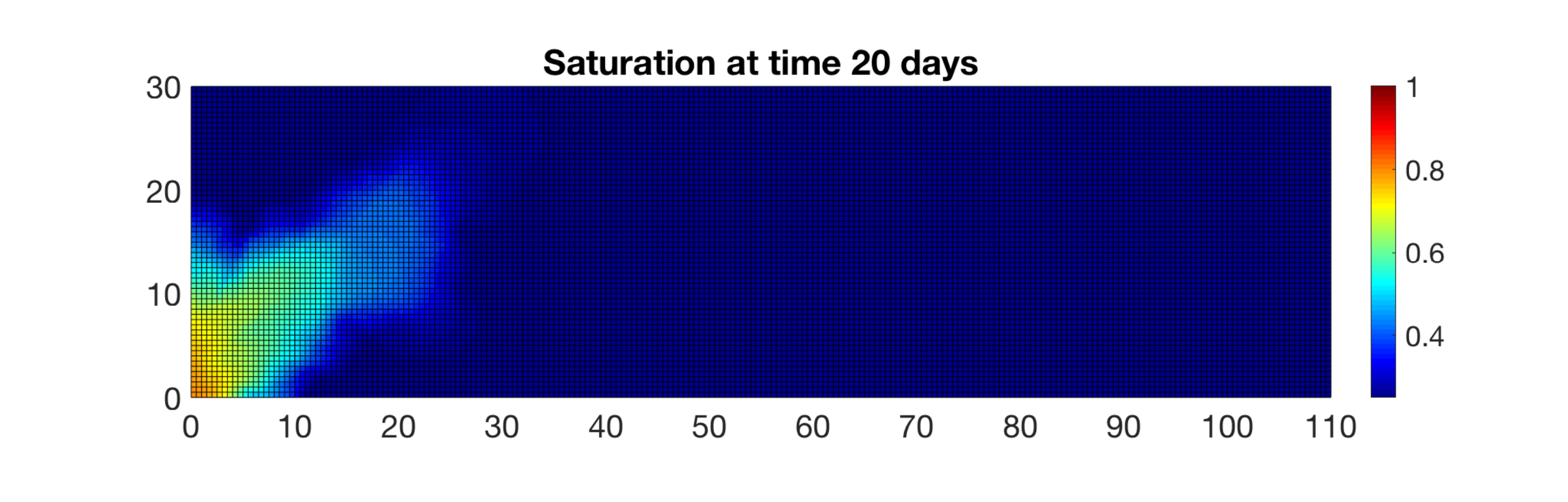}
\includegraphics[width=7.5cm,trim=2.5cm 1cm 2.5cm 0.5cm, clip]{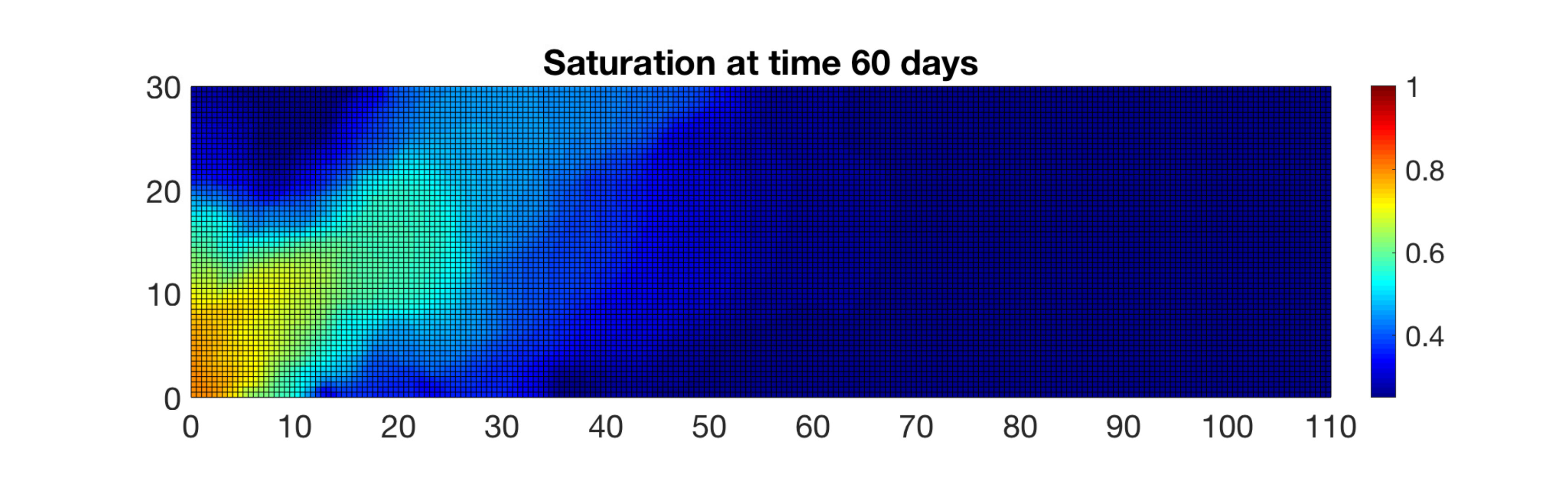}
\caption{Saturation distribution (top) and domain decomposition map (middle) for the adaptive space-time approach, and fine scale simulation comparison at 20 and 60 days (layer 20).}
\label{fig:l202phsat}
\end{center}
\end{figure}

We observed a computational speedup of approximately 25 and 40 times with our dynamic mesh refinement in space and time approach when compared to fine scale simulations for layer 20 (Gaussian permeability distribution) and 37 (channelized permeability distribution) cases, respectively. We attribute this difference in speedup to the variation velocities at the saturation front ($S_{w}^{*}$) for these two cases. Let us assume, for the sake of argument, that the saturation front can be characterized by an average velocity. Then along this saturation front, the velocities at different spatial and temporal locations vary about this average. For the Gaussian permeability distribution, this variation is not too large and therefore the front is mostly classified as identifier 1 (fine in space and time) which is computationally intensive. On the other hand, for the channelized permeability distribution this variation is large due to high contrast in permeability values and therefore a narrower region is classified as identifier 1 resulting in lower computational costs. Further for the channelized case, the saturation front sweeps a smaller area and consequently a large part of the domain is classified as identifier 4 (coarse in space and time) and thus a lower computational cost compared to Gaussian distribution. The identifier maps in Figures \ref{fig:l372phsat} and \ref{fig:l202phsat} shows a visual representation of this argument at 20 and 60 days.

\section{Conclusions}
An adaptive space-time, domain decomposition approach is presented for balancing computational loads associated with solving multiphase, flow and transport problems in porous medium. We also described an efficient space-time monolithic solver for this approach that does not require subdomain iterations. A standalone, serial prototype was developed for testing and preliminary computational benchmarking purposes. The explicit error estimators (normalized non-linear residuals) were used to indicate problem areas that pose convergence issues for the non-linear solver. Our adaptive approach then relies upon a combination of explicit error estimators and delta change in the quantity of interest in space and time to divide the reservoir into subdomains with different spatial and temporal mesh refinements. This allows us to circumvent non-linear solver convergence issues towards promoting computational efficiency without requiring time step size reduction for the entire domain in the event of a convergence failure. A rigorous derivation of problem specific explicit and a-posteriori error estimators equipped with appropriate norms is reserved for a future work.  We observed a computational speedup of 25 to 40 times when our framework was applied to a slightly compressible, two-phase flow problem in a heterogeneous porous media compared to conventional fine-scale simulations. The time-concurrent, solution algorithm reduces the serial nature of the conventional, time-marching solution algorithms. This renders a massively parallel framework where parallel in time, linear solvers and preconditioners can be used without compromising computational efficiency in the near future.

\section{Acknowledgements}
This work was supported by Department of Energy (DOE), Center for Frontiers for Subsurface Energy Security (CFSES) grant DE-SC000111 and National Science Foundation (NSF), BIGDATA: Collaborative Research Award 1546553. We would also like to thank Center for Subsurface Modeling (CSM) industrial affiliates for their continued support. 

\bibliography{ref}
\bibliographystyle{plain}

\end{document}